\documentclass[a4paper,11pt]{article}

\usepackage{amsmath}
\usepackage{amsfonts}
\usepackage{amssymb}
\usepackage{amsthm}
\usepackage{graphicx}
\usepackage{verbatim}
\usepackage{authblk}

\newtheorem{thm}{Theorem}[section]

\newtheorem{prop}[thm]{Proposition}
\theoremstyle{definition}
\newtheorem{defin}[thm]{Definition}

\newcommand{\eps}{\varepsilon}

\newcommand{\N}{\mathbb{N}}
\newcommand{\Q}{\mathbb{Q}}
\newcommand{\R}{\mathbb{R}}
\newcommand{\T}{\mathbb{T}}

\newcommand{\Graph}{\text{Graph}}

\author[1]{Angel Jorba}
\author[2]{Pau Rabassa}
\author[1]{Joan Carles Tatjer}
\affil[1]{\mbox{Departament of Matem\`atica Aplicada i An\`alisi,}
\mbox{Universitat de Barcelona, Barcelona, Spain}\vspace{2mm}}
\affil[2] {\mbox{Johann Bernoulli Institute for Mathematics and Computer Science,}
\mbox{University of Groningen, Groningen, The Netherlands}}

\date{}

\title{Period doubling and reducibility in the quasi-periodically 
forced logistic map\thanks{This work has been supported by the MEC grant MTM2009-09723
and the CIRIT grant 2009 SGR 67.}
}
\begin{document}
\maketitle

\begin{abstract}
We study the dynamics of the Forced Logistic Map
 in the cylinder. We compute a bifurcation diagram in terms 
of the dynamics of the attracting set. 
Different properties of the attracting set
are considered, as the Lyapunov exponent  and, in the case of having
a periodic invariant curve, its period and its reducibility. 
This reveals that the parameter values for which the invariant
curve doubles its period are contained in regions of the parameter
space where the invariant curve is reducible. Then we present 
two additional studies to explain this fact. In first place we consider the
images and the preimages of the critical set (the set where the
derivative of the map w.r.t the non-periodic coordinate is equal to zero).
Studying these sets we construct constrains in the parameter 
space for the reducibility of the invariant curve. In second place 
we consider the reducibility loss of the invariant curve as
codimension one bifurcation and we study its interaction with
the period doubling bifurcation. This reveals that, if the 
reducibility loss and the period doubling bifurcation curves meet, 
they do it in a tangent way. 
\end{abstract}

\tableofcontents

%%%%%%%%%%%%%%%%%%%%%%%%%%%%%%%%%%%%%%%%%%%%%%%%%%%%%%%%%%%%%%%%%%%%%%%%%

%%%%%%%%%%%%%%%%%%%%%%%55555
%%%%%%%%%%%%%%%%%%%%%%%%
\section{Introduction}
\label{section}
%%%%%%%%%%%%%%%%%%%%%%%%
%%%%%%%%%%%%%%%%%%%%%%%%

We  focus on the study of the quasi-periodically
(q.p. for short) Forced Logistic Map (FLM for short). The FLM is a two
parametric map in the cylinder 
%$(\T \times \R)$ 
where the dynamics in the periodic component is a rigid rotation 
and the dynamics in the other component
is the logistic map plus a quasi-periodic forcing term. 
This map appears in the literature in different
contexts, usually related with the destruction of invariant curves.
For example, in \cite{HH94} it was introduced as an example
where SNAs (Strange Non-chaotic Attractors) were created through
a collision between stable and unstable invariant curves. Since then,
different routes for the destruction of invariant
curves have been explored for this map, for instance see \cite{PMR98}
and references therein. Some other recent studies are 
\cite{AKIL07,AC09, Bje09,FH10, Jag03,JT05}.

On the other hand, the FLM is also related to the truncation of
period doubling cascades. It is well known that the one dimensional 
logistic map exhibits an infinite cascade of period doubling bifurcations
which leads to chaotic behavior. Moreover this infinite cascade extends 
to a wider class of unimodal maps.  
But when some q.p. forcing is added, the number of period doubling 
bifurcations of the invariant curves is finite. 
This phenomenon of finite period doubling cascade has been observed in different
applied and theoretical contexts. In the applied context it has been
observed in a truncation of the Navier-Stokes flow \cite{Fran83,Veen05} or in
a periodically driven low order atmosphere model \cite{BSV02}. In the
theoretical context, it has also been reported in different maps
which were somehow built to have period doubling cascades \cite{ACS83,Kan83}, and more
recently in the analysis of the Hopf-saddle-node bifurcation \cite{BSV09}.
Actually, in \cite{Kan83} the FLM itself is given as a model for the
truncation of the period doubling bifurcation cascade.

The study presented here is more concerned with
the mechanisms which cause the truncation of the period doubling
bifurcation cascades than with the possible existence
of SNAs for this family of maps. Concretely, we show (numerically) 
that the reducibility has the role of confining the period doubling 
bifurcation in closed regions of the parameter space. In the remainder 
of this article we focus on the shape of this reducibility regions. 
In \cite{JRT11a, JRT11b, JRT11c} we will use the reducibility loss bifurcation to 
study the self renormalizable properties of the bifurcation 
diagram and how the Feigenbaum-Collet-Tresser renormalization theory 
can be extended to understand it. See also \cite{Rab10} for a united 
exposition of the present paper with the other three cited before.

%\fbox{Review the summary below} 

This paper is structured as follows. In section 
\ref{section invariant curves in q.p systems} 
we review some concepts and results concerning
the continuation of invariant curves for a quasi-periodic 
forced maps. We also look at the concrete case
when the map is uncoupled. 

In section \ref{chapter forcel logistic map} we focus on the 
 dynamics of the FLM. First we review some computations
which can be found in the literature. Then, we do a study 
of the parameter space in terms of the dynamics of the attracting set of the map.
For this study different properties of the attracting set
are considered, as the value of the Lyapunov exponent  and, in the case of having
a periodic invariant curve, the period. Differently to other works, in
our study the reducibility of the invariant curves
has been also taken into account. This reveals interesting information,
for example we observe that the parameter values for which the invariant
curve doubles its period is contained in regions of the parameter
space where the invariant curve is reducible. The subsequent sections are
developed with the aim of understanding the results presented
in this section.

In section \ref{section obstruction to reducibility} we consider the 
images and the preimages of the critical set (this set is the set where the 
derivative of the map w.r.t the non-periodic coordinate is equal to zero). We
also consider the continuation in the parameter space of the invariant curve
which comes from one of the fixed points of the logistic map. Doing a
study of the preimages of the critical set we construct forbidden regions
in the parameter space for the reducibility of the invariant curve. 
In other words, we give some constrains on the reducibility of the invariant 
curve.

In section \ref{section period doubling and reducibility}
the reducibility loss of the invariant curve is considered as a
codimension one bifurcation and then we study its interaction with
the period doubling bifurcation. The study done here is not particular
for the FLM. 
In this section we also give a general model for the
reducibility regions enclosing the period doubling bifurcation
observed in the parameter space of section \ref{section parameter 
space and reducibility}.
%This study is not done in the FLM but from a more general point of view. 

In section \ref{section summary and conclusions}
 we summarize the results obtained in the previous ones.
We analyze  again the bifurcations diagram obtained in section 
\ref{section parameter space and reducibility} taking into account 
the results obtained in sections  \ref{section obstruction to reducibility} and 
\ref{section period doubling and reducibility}. 

%%%%%%%%%%%%%%%%%%%%%%%%%%%%%%%%%%%%%%%%%%%%%%%%%%%%%%%%%%%%%%%%%%%%%%%555
%%%%%%%%%%%%%%%%%%%%%%%%%%%%%%%%%%%%%%%%%%%%%%%%%%%%%%%%%%%%%%%%%%%%%%%555
\section{Invariant curves in quasi-periodically forced systems}
\label{section invariant curves in q.p systems}
%%%%%%%%%%%%%%%%%%%%%%%%%%%%%%%%%%%%%%%%%%%%%%%%%%%%%%%%%%%%%%%%%%%%%%%555
%%%%%%%%%%%%%%%%%%%%%%%%%%%%%%%%%%%%%%%%%%%%%%%%%%%%%%%%%%%%%%%%%%%%%%%555

In this section we briefly review some of the key definitions 
and results on the theory of invariant curves
 in quasi-periodically forced maps. These definitions will 
be useful for the forthcoming analysis of the 
dynamics of the FLM. 

\subsection{Basic definitions}

A {\bf quasi periodically forced one dimensional map}  is a map of the form
\begin{equation}
\label{q.p. forced map general}
\begin{array}{rccc}
F: & \T \times \R & \rightarrow & \T \times \R  \\
& \left( \begin{array}{c} \theta \\ x \end{array} \right) &
\mapsto &  \left( \begin{array}{c}  \theta + \omega \\ f(\theta,x)
\end{array} \right)
\end{array}
\end{equation}
where $f\in C^r(\T\times \R ,\R)$ with $r\geq 1$ 
%(including $r=\infty$ and the analytic cases) 
and the parameter $\omega\in \T\setminus\Q$. 
%For some results we will need additional requirements on $\omega$ and 
%on the differentiability of $f$. 
%Note that any map like (\ref{q.p. forced map general}) can be identified with an element 
%of the space $C^r(\T\times \R ,\R)\times \T$. 

Given a quasi-periodically  forced map as above, we have that it determines 
a dynamical system in the cylinder, explicitly defined
as 
\begin{equation}
\left.
\begin{array}{lr}
\bar{\theta}=\theta + \omega, \\
\bar{x}= f(\theta,x), 
\end{array}
\right\}
\label{q.p. forced system}
\end{equation}

%Given a quasi periodically forced system 
%like (\ref{q.p. forced system}) we are interested in the invariant sets
%it may have. Since we are assuming $\omega$ to be irrational 
%the simplest invariant sets are invariant curves. 

\begin{defin}
Given a continuous  function  $u:\T \rightarrow \R$ we will say that $u$ is an 
{\bf invariant curve}  of  (\ref{q.p. forced system}) if, and only if,
\begin{equation}
u(\theta+ \omega) = f(\theta, u(\theta)), \quad \forall \theta \in \T.
\label{invariance equation}
\end{equation}

The value $\omega$ is known as the {\bf rotation number} of $u$.
\end{defin}

An equivalent way to define invariant curve, is to require  the set 
$\{(\theta,x)\in \T\times\R |\thinspace x=u(\theta)\}$ to be invariant 
by $F$, where $F$ is the function defined by (\ref{q.p. forced map general}). 

On the other hand, if we consider the map $F^n$ we have that it  is also 
a quasi-periodically 
 forced map. Given a function  $u:\T \rightarrow \R$, we will 
say that $u$ is a $n$-periodic invariant curve of $F$ if the set 
$\{(\theta,x)\in \T\times\R |\thinspace x=u(\theta)\}$ is invariant by $F^n$ 
(and there is no smaller $n$ satisfying such  condition). 

Since a periodic invariant curve of a map $F$ is indeed 
an invariant curve of $F^n$, any result for 
invariant curves can be extended to periodic 
invariant curves. %On the other hand, let us remark that 
%if we have $f\in C^r(\T\times \R,\R)$, and $u \in C^0(\T,\R)$ an
%attracting invariant curve of (\ref{q.p. forced system}), then
%one can prove (see \cite{HL05}) that $u$ is indeed in $C^r(\T,\R)$.

Given $x=u_0(\theta)$ an invariant curve of (\ref{q.p. forced system}), 
its linearized normal behavior is described by the following linear 
skew product:
\begin{equation}
\label{linear skew product}
\left.
\begin{array}{rcl}
\bar{\theta}& = &\theta  + \omega, \\
\bar{x} & = & a(\theta) x, 
\end{array}
\right\}
\end{equation}
where $a(\theta) = D_x f( \theta, u_0(\theta))$ is also of class $C^r$,
$x \in \R$ and $\theta \in \T$. We will assume that the invariant curve is 
not degenerate, in the sense that the function $a(\theta)$ is not identically
zero.

\begin{defin}
\label{definition reducibility}
The system (\ref{linear skew product}) is called {\bf reducible} if, and only if, 
there exists a change of variable $x=c(\theta) y$, 
continuous with respect to $\theta$, such that (\ref{linear skew product}) 
becomes
\begin{equation}
\label{reduced linear skew product}
\left.
\begin{array}{rcl}
\bar{\theta}& = &\theta  + \omega, \\
\bar{y} &= & b y, 
\end{array}
\right\}
\end{equation}
where $b$ does not depend on $\theta$. The constant $b$ is 
called the {\bf multiplier} of the reduced system.
\end{defin}

In the case that $a(\cdot)$ is a $C^\infty$ function and $\omega$ 
is Diophantine (see Proposition 1 in \cite{JT05}), 
the skew product (\ref{linear skew product})
 is reducible  if, and only if,
 $a(\cdot)$ has no zeros, see Corollary 1 of \cite{JT05}. 
Actually, % under the same hypothesis  
the reducibility loss can be characterized 
as a codimension one bifurcation. 
\begin{defin}
\label{definition of reducibility loss bifurcation}
Let us consider a one-parametric family of linear skew-products
\begin{equation}
\label{1-d family of skew products}
\left.
\begin{array}{rcl}
\bar{\theta}&= & \theta  + \omega, \\
\bar{x}& = & a(\theta, \mu) x, 
\end{array}
\right\}
\end{equation}
where $\omega$ is Diophantine and $\mu$ belongs to an open 
set of $\R$ and $a$ is a $C^\infty$ function of $\theta$ and $\mu$. We 
will say that the system (\ref{1-d family of skew products})
undergoes a {\bf reducibility loss bifurcation} at $\mu_0$ if 
\begin{enumerate}
\item $a(\cdot, \mu)$ has no zeros for $\mu<\mu_0$,
\item $a(\cdot, \mu)$ has a double zero at $\theta_0$ for $\mu=\mu_0$,
\item $\frac{d}{d\mu} a (\theta_0, \mu_0) \neq 0$. 
\end{enumerate}
\end{defin}

On the other hand, consider a system like (\ref{q.p. forced system}) with 
$f$ a $C^\infty$ function, which depends (smoothly) on a one dimensional 
parameter ($f=f_\mu$). 
Assume also that we have an invariant curve $u=u_\mu$ of the system. 
We will say that the invariant curve undergoes a reducibility 
loss bifurcation if the system (\ref{linear skew product}) associated to the invariant curve
($a(\theta) = a(\theta,\mu) = D_x f_\mu(\theta, u_\mu(\theta))$)
undergoes a reducibility loss bifurcation as a system of linear skew-products.

Given a map like (\ref{linear skew product}) 
we have that, due to the rigid rotation in
the periodic component,  one of Lyapunov exponents is equal to zero
(see \cite{AC09}).  Then the definition of the Lyapunov exponent 
can be suited to the case of linear skew-products as follows.

\begin{defin}
\label{definitio Lyapunov exponent skew-products}
If $\theta \in \T$, we define the {\bf Lyapunov exponent} of (\ref{linear
skew product}) at $\theta$ as 
\begin{equation}
\label{Lyapunov exponent limit}
\lambda(\theta) = \limsup_{n\rightarrow \infty} \frac{1}{n} 
\ln \left| \prod_{j=0}^{n-1} a(\theta  + j \omega) \right|.
\end{equation}
We also define the {\bf Lyapunov exponent of the skew product} (\ref{linear
skew product}) as
\begin{equation}
\label{Lyapunov exponent integral}
\Lambda =  \int_{0}^{1} \ln|a(\theta)|d\theta.
\end{equation}
\end{defin}

If $\Lambda$ is finite then, applying the Birkhoff Ergodic Theorem 
%(\cite{KH96}), 
we have that 
the $\limsup$ in (\ref{Lyapunov exponent limit}) is in fact a limit and  
$\lambda(\theta) = \Lambda$
for Lebesgue a.e. $\theta \in \T$. 
If $a(\theta)$ never vanishes, 
the $\limsup$ in (\ref{Lyapunov exponent limit}) is again a limit and
coincides with $\Lambda$, but now for all $\theta \in \T$.

Now, consider  a map like (\ref{q.p. forced map general}). If there exists
an invariant compact cylinder where $f$ is 
monotone and has a negative 
Schwarzian derivative with respect to $x$, then 
J\"ager has proved the existence of invariant curves, 
see \cite{Jag03} for details. On the other hand in \cite{JT05} 
a result on the persistence of invariant 
curves is given, in terms of the reducibility and the Lyapunov exponent
 of the curve.

%%%-----------------------------------------------%%% 
\subsection{Quasi-periodically forced maps which are uncoupled}
%%%-----------------------------------------------%%% 

In this subsection we turn our attention to 
the maps which are of the same 
class of the FLM, in the sense that the quasi-periodic function can 
be written as a one dimensional function plus a quasi-periodic term.

\begin{defin}
Given a map like (\ref{q.p. forced map general}) we will say the the map 
{\bf is uncoupled} if $f$ does not depend on $\theta$, i.e. $f(\theta,x)=f(x)$
\end{defin}
Note that if a map $F$ is uncoupled, the $F^n$ also does. 

%A first property of these maps is the following
\begin{prop}
\label{prop persistence fixed points}
Let $F_\mu$ be a one dimensional family of maps 
like (\ref{q.p. forced map general}) such 
that for a fixed value $\mu_0$ the map $F_{\mu_0}$  is uncoupled, that is
$F_{\mu_0}(\theta,x) = (\theta+\omega, f_{\mu_0}(x))$. 
Then any hyperbolic fixed point of $f_{\mu_0}$
extends to an invariant curve of the system
for $\mu$ close to $\mu_0$.
\end{prop}

\begin{proof}
Suppose that there exists a hyperbolic fixed point $x_0\in \R$ of 
$F_{\mu_0}$. We have that it can be seen as an invariant curve 
$u_0:\T \rightarrow \R$ of $F_{\mu_0}$, with $u_0(\theta) = x_0$ 
for any $\theta\in\T$. The skew product (\ref{linear skew product}) 
associated  to the invariant curve has as a multiplier 
$a(\theta) = \frac{\partial}{\partial x}f(\theta,x_0) = f_{\mu_0}'(x_0)$, 
which actually does not depend on $\theta$. Concretely we have 
that the system is reducible. 
  
Now we can apply the theory exposed in section 3.3 of \cite{JT05}. 
Assume first that $f_{\mu_0}'(x_0) \neq 0$. 
We have that a curve is persistent by perturbation if $1$ does not 
belong to the spectrum of the transfer operator associated to the curve. Since the curve is 
reducible we have that the spectrum is a circle of modulus $|a|$. 
Using that the fixed point $x_0$ is hyperbolic we have that $|a|\neq 1$, 
therefore $1$ does not belongs to the spectrum. When $f_{\mu_0}'(x_0) =0$
we have that the spectrum collapses to $0$. Then $1$ does not 
belong to the spectrum of the transfer operator either. 
\end{proof}

Note that, by considering $F^n$ for the periodic case, the result 
extends to any periodic point of the uncoupled system. 

This last proposition can be also proved using the normal hyperbolicity 
theory \cite{HPS70}. But this theory is only valid for diffeomorphisms, 
then  the case of $ f_{\mu_0}'(x_0) = 0$ is not
included.

%%%%%%%%%%%%%%%%%%%%%%%%%%%%%%%%%%%%%%%%%%%%%%%%%%%%%%%%%%%%%%%%%%%%%%%%%

%%%%%%%%%%%%%%%%%%%%%%%%%%%%%%%%%%%%%
\section{The Forced Logistic Map}
\label{chapter forcel logistic map}
%%%%%%%%%%%%%%%%%%%%%%%%%%%%%%%%%%%%%%

The FLM is a map in the cylinder $\T\times \R$
defined as
\begin{equation}
\label{FLM}
\left.
\begin{array}{rclcl}
\bar{\theta} & = & r_\omega(\theta) & = & \theta + \omega  ,\\
\bar{x} & = & f_{\alpha,\eps} (\theta,x) & = & \alpha x(1-x)(1+ \eps \cos(2\pi \theta)) ,
\end{array}
\right\}
\end{equation}
where $(\alpha,\eps)$ are parameters and $\omega$ a fixed Diophantine
number (typically in our study it will be the golden mean).

%%%%%%%%%%%%%%%%%%%%%%%%%%%%%%%%%%%%%%%%%5
\subsection{Basic study of the dynamics}
\label{section basic study of FLM}
%%%%%%%%%%%%%%%%%%%%%%%%%%%%%%%%%%%%%%%%%%%5

Note that the FLM (\ref{FLM}) is a q.p. forced
system like (\ref{q.p. forced system}),
which depends on two parameters $\alpha$ and $\eps$. Moreover,
we have that the function $f_{\alpha,\eps}$ which defines the map
can be written as a logistic map plus a q.p. forcing term.
In other words, we have that
\[ f_{\alpha,\eps}(\theta,x)= \ell_\alpha(x) + h_{\eps,\alpha}(\theta,x), \]
 where $\ell_\alpha(x)=\alpha x(1-x)$ (the logistic map) and
$h_{\eps,\alpha}(\theta,x) = \eps \alpha x(1-x)\cos(2\pi\theta)$ 
(which is zero when $\eps$ is). Then proposition 
\ref{prop persistence fixed points} is applicable to
the map. In some cases it will be convenient to work in a compact
domain. In this case note that when $0\leq \alpha(1+|\eps|) \leq 4$ the
compact cylinder $\T \times[0,1]$ is invariant by the map.

The one dimensional logistic map is known to exhibit a period
doubling bifurcation cascade, that is a sequence of infinitely 
many period doubling bifurcations accumulating to a concrete parameter value. 
Applying proposition
\ref{prop persistence fixed points}, we have that the period
doubling cascade persists in the sense that,
fixed a hyperbolic periodic orbit of period  $2^n$ of the logistic map,
there exists a sufficiently small $\eps$  such that the
prescribed periodic orbit extends to a periodic invariant curve for 
$\eps\leq \eps_0$.

\begin{figure}[t]
\begin{center}
\includegraphics[width=5cm]{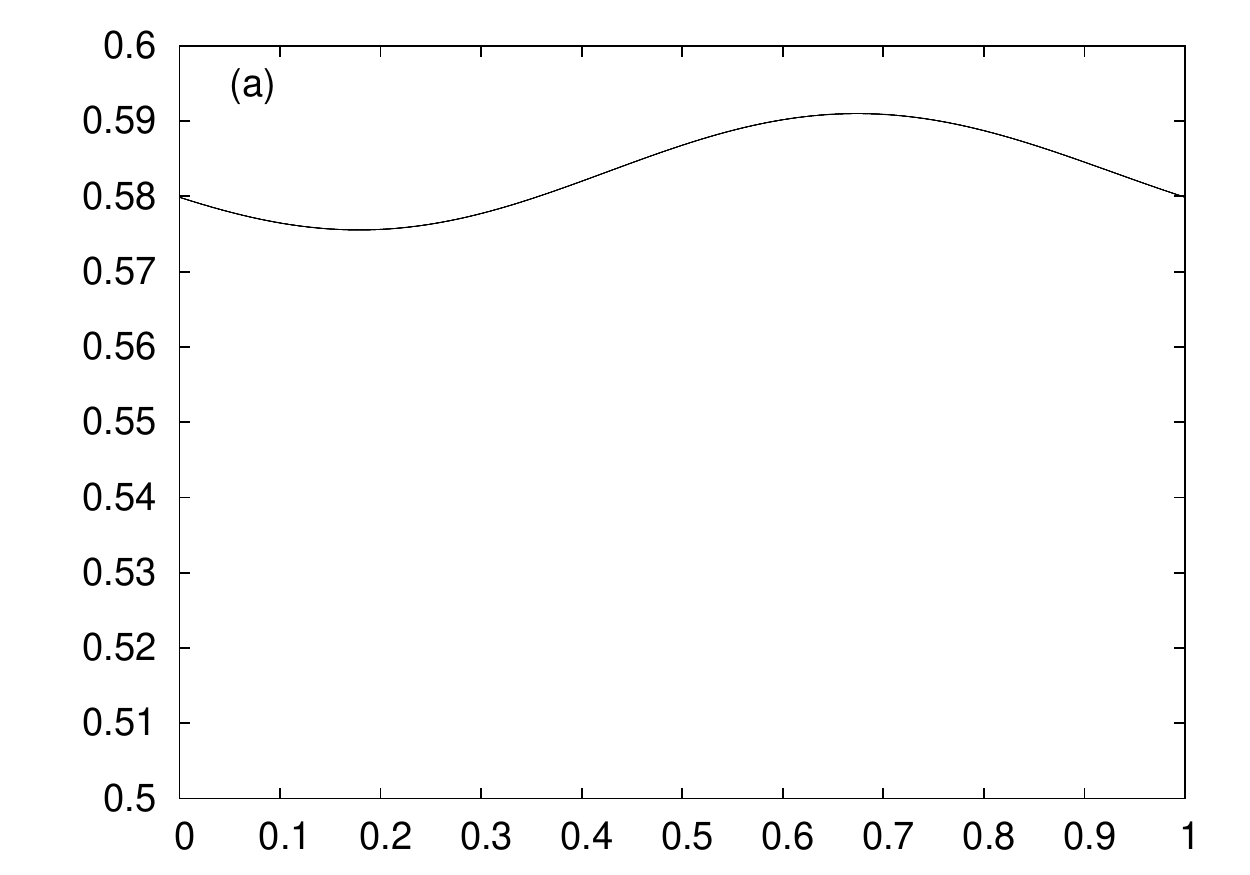}
\includegraphics[width=5cm]{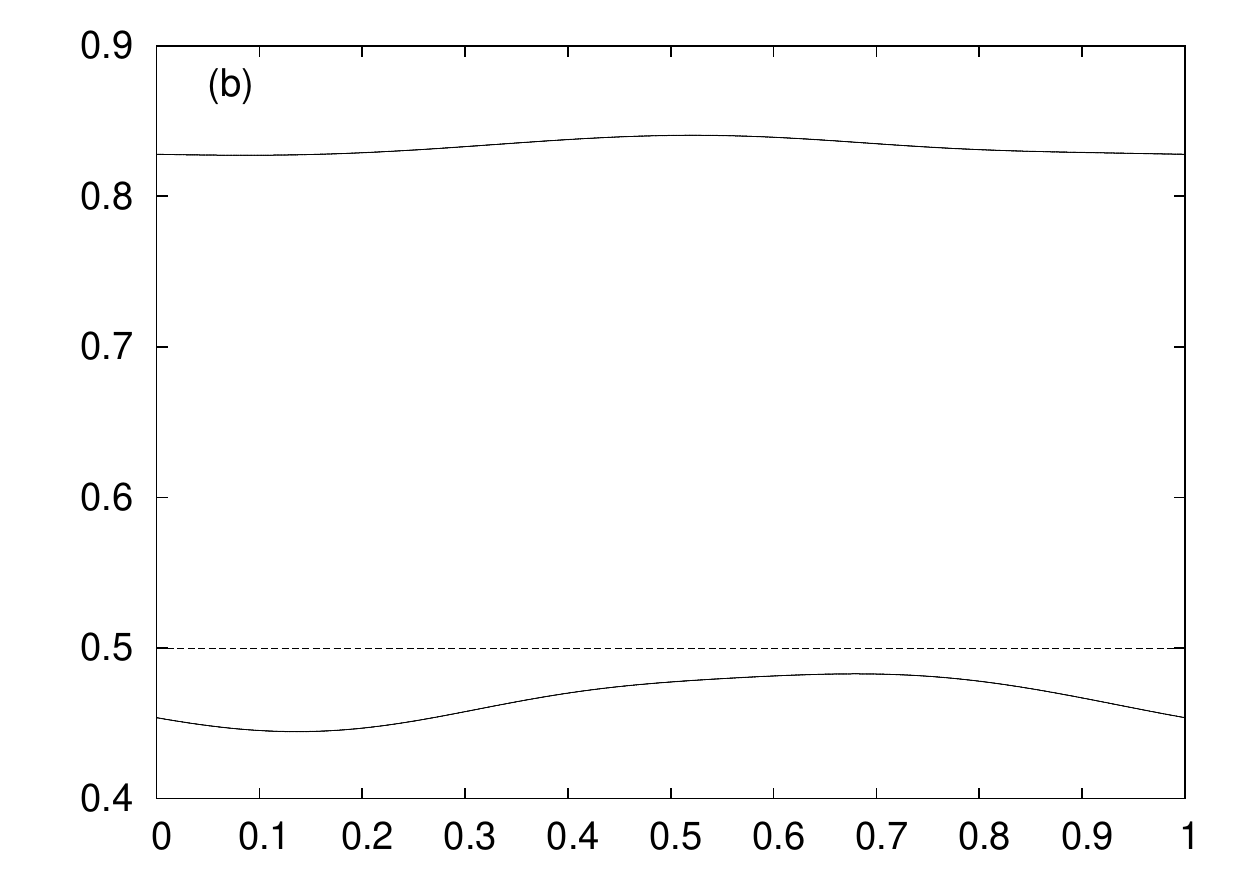}
\includegraphics[width=5cm]{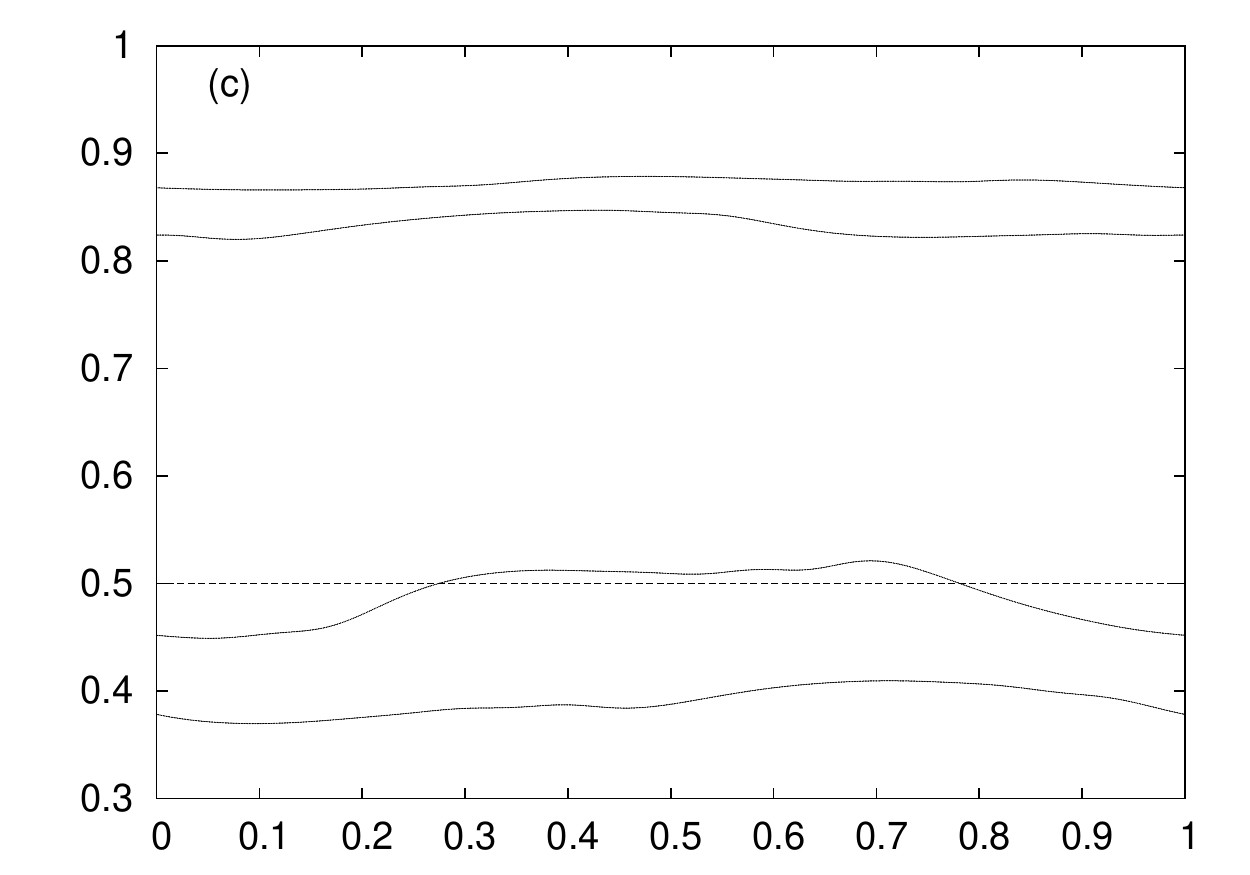}
\includegraphics[width=5cm]{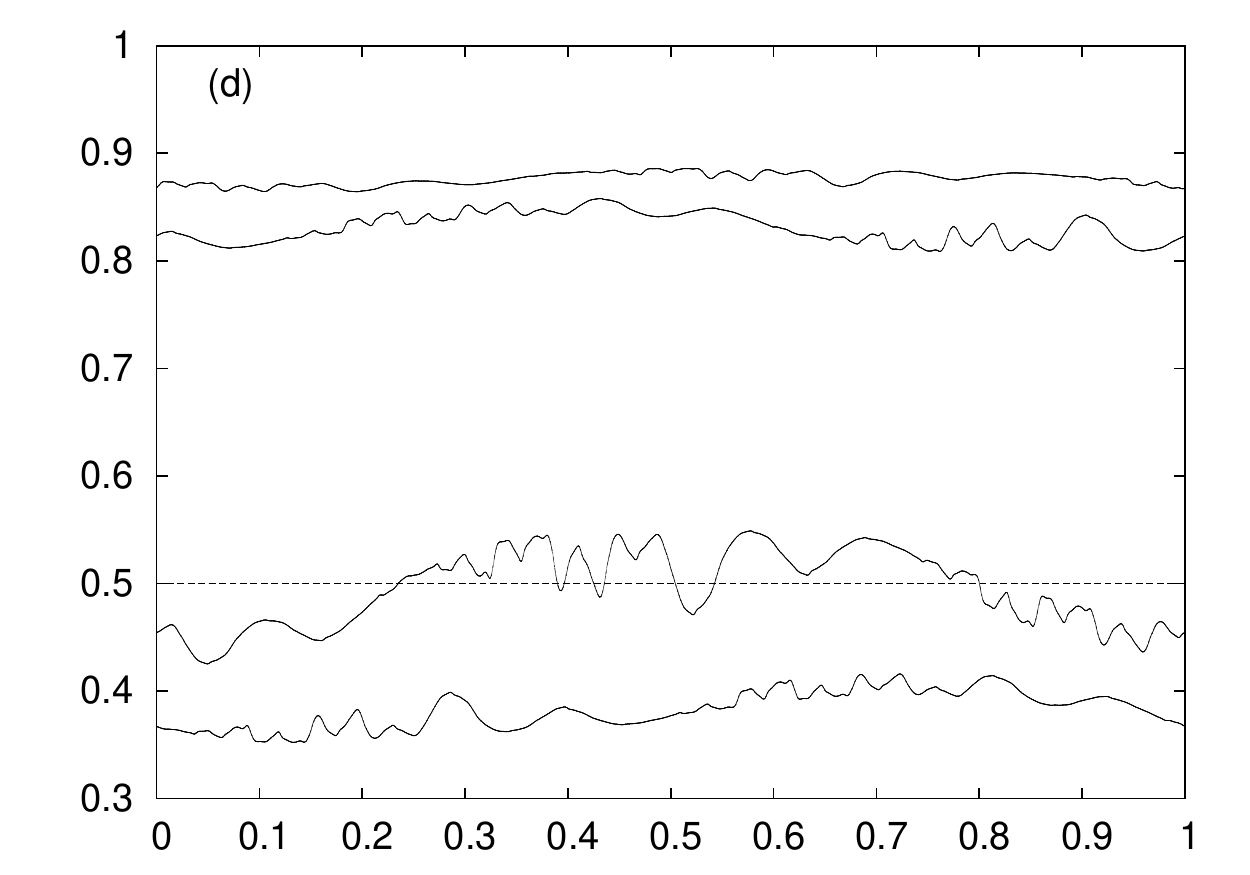}
\includegraphics[width=5cm]{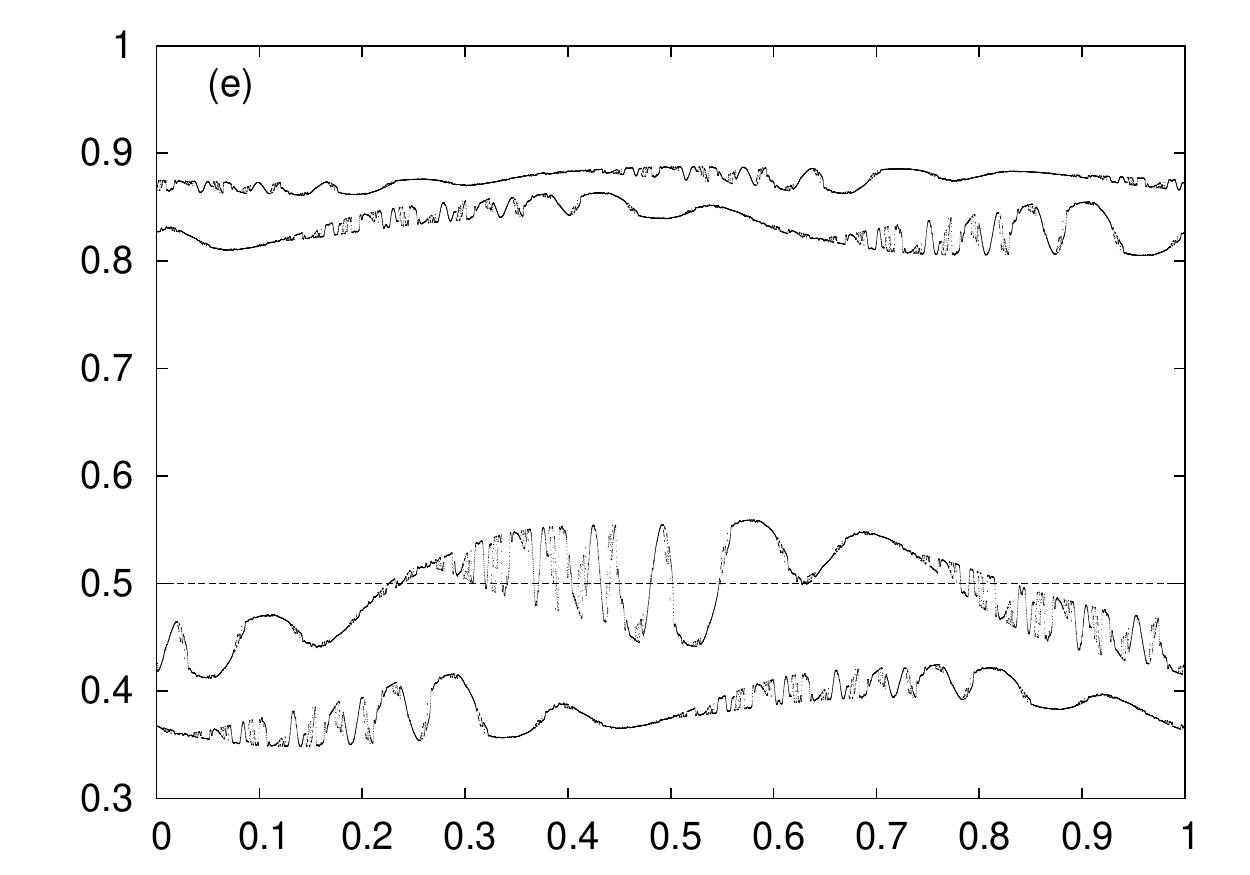}
\includegraphics[width=5cm]{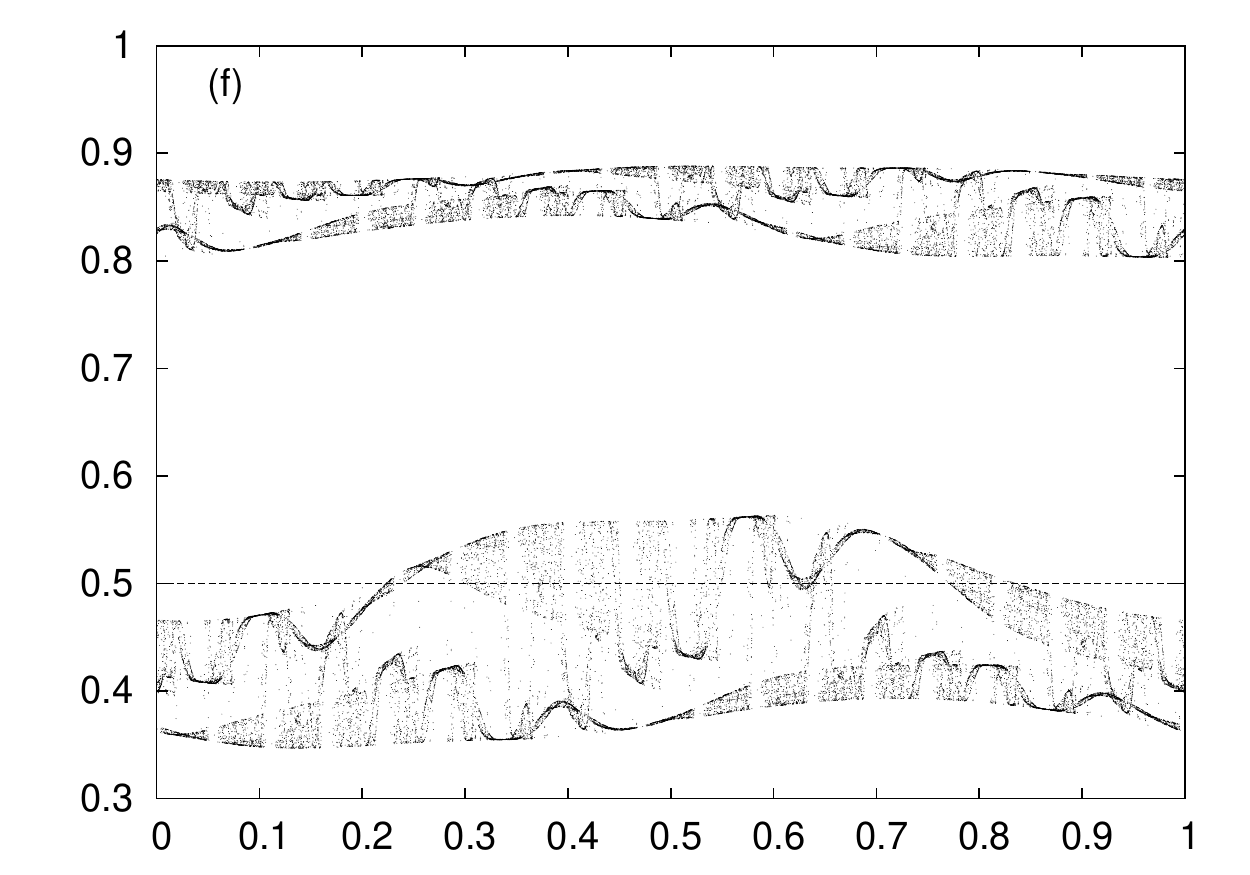}
\caption{The attractor of the FLM (\ref{FLM}),
for different parameters of $\alpha$ when $\eps$ is $0.01$.
The vertical axis corresponds to the $x$ variable and the horizontal 
one to $\theta$. 
The values of $\alpha$ (reading from left to right and top to bottom)
are $2.4$, $3.35$, $3.5$, $ 3.522$, $3.527$ and $3.529$.
The values of the Lyapunov exponent for each attractor are plot
in figure \ref{Lyapunov exponent eps fix}.
\label{exemples fractalitzacio a periode 1}}
\end{center}
\end{figure}

\begin{figure}[t]
\begin{center}
\includegraphics[width=7.5cm]{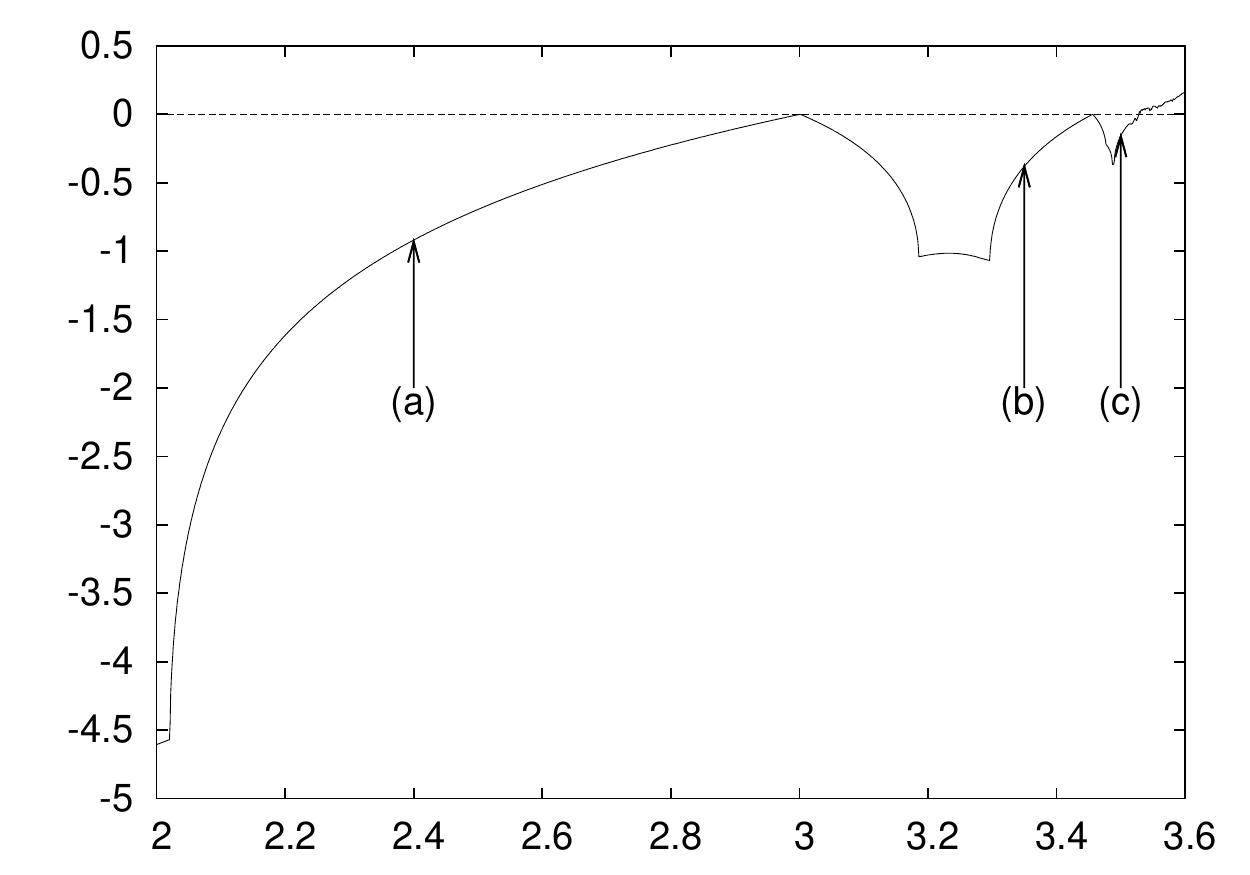}
\includegraphics[width=7.5cm]{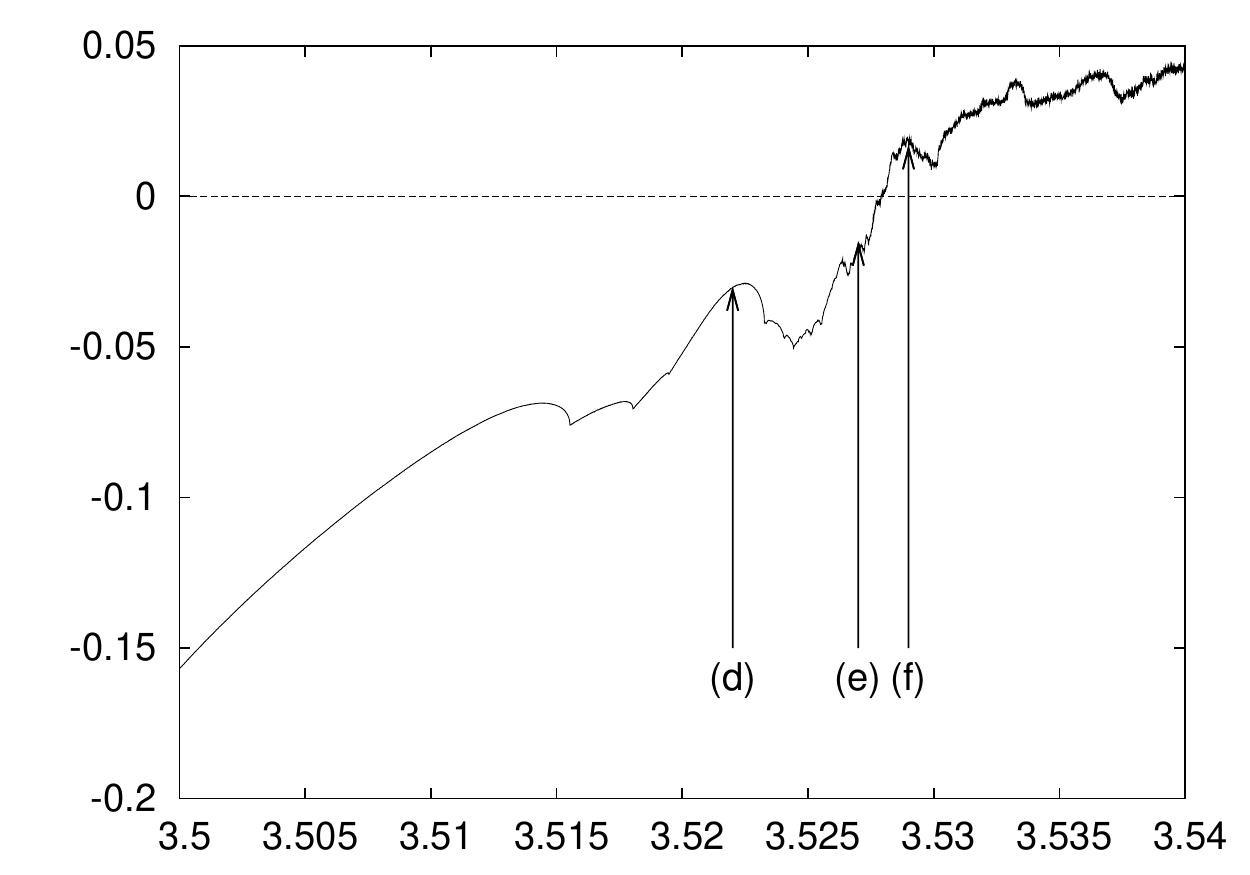}
\caption{
Graph of the Lyapunov exponent (vertical axis) of  the attracting set of
the FLM (\ref{FLM}) as a function of $\alpha$ (horizontal axis), 
for $\eps= 0.01$. The values corresponding to
the attractors displayed in figure
\ref{exemples fractalitzacio a periode 1} are also shown.}
\label{Lyapunov exponent eps fix}
\end{center}
\end{figure}

A natural question to ask is what happens in the converse sense. That is,
consider the parameter $\eps$ fixed (equal to a small value) and let
the parameter $\alpha$ increase. What
happens to the period doubling bifurcation cascade of the
logistic map? The answer to
this question can be easily guessed with some simple
computations.

In figure \ref{exemples fractalitzacio a periode 1} we have plotted
the attractor of the FLM for several values of $\alpha$ when $\eps$
is $1/100$. We have also computed  the Lyapunov exponent
of the attractor,  as a
function of $\alpha$ in figure \ref{Lyapunov exponent eps fix}.
The  attractor has been computed by forward iteration of the map,
using a transient and then plotting a certain number
of iterates to obtain the approximation of the attractor. For the
computation of the Lyapunov exponent we have used
its definition as a limit (see (\ref{Lyapunov exponent limit})). 
%\ref{subsection Lyapunov exponents}

For the first three values of $\alpha$ we observe how
the attracting invariant curve doubles its period twice. After
two period doublings the curve begins to get more and
more wrinkled until a strange\footnote{
We will say that a set is strange if it is neither a
differentiable curve nor a finite union of them.} attractor seems to appear.
This scheme is known as the truncation of the period doubling bifurcation and
it is reproduced for any (arbitrarily small but fixed) value of $\eps$.

%fbox{DIR ALGO MES PER LLIGAR-HO AMB TOT LO QUE VE!} 

In figure \ref{Lyapunov exponent eps fix} we see that in the graph of the
Lyapunov each period doubling of the attracting
curve corresponds to a value of $\alpha$ where the exponent
becomes zero. Note that between
one period doubling and the next one  the Lyapunov exponent of the
system (as a function of $\alpha$)  has some critical values where
its derivative goes to minus infinity. This behavior of the Lyapunov
exponent is due to a loss and a recovery of reducibility of
the attracting curve (see \cite{JT05}).
If one increases $\alpha$,
after a certain number of period doublings, there are no more bifurcations of
this kind. Then, the attracting curve wrinkles until it becomes a strange 
set. In the graph of the Lyapunov exponent in figure \ref{Lyapunov exponent eps fix} 
we observe that it goes to zero,
but now it does it in a sharp way and then it  becomes
positive. This is a typical behavior also reported in other maps where
there is a truncation of the period doubling sequence \cite{Kan83,Veen05,BSV09}.
Therefore the FLM has a lot interest as a simple model for this process.

The destruction of the invariant curve shown above
 is indeed a subtle problem. For certain parameter values
%during  this process but just before the Lyapunov exponent getting zero, 
the numerical computations of the attractor produce a set which seems  strange. One
might think that it is a SNA, but doing computations with
higher accuracy, after several magnifications, the numerical
approximation of the invariant curve does not seem strange anymore,
but is an incredibly wrinkled invariant
curve (see \cite{HS05,JT05}). Actually, for certain parameter values
one is not able to conclude (numerically) whether the curve is strange or not.

\begin{figure}[t]
\begin{center}
\includegraphics[width=7.5cm]{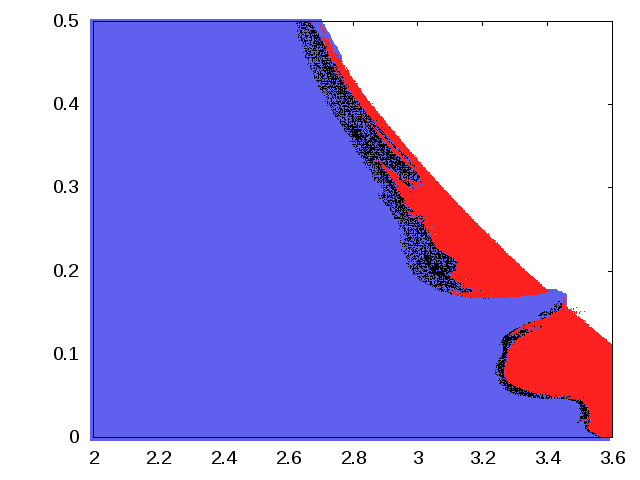}
\includegraphics[width=7.5cm]{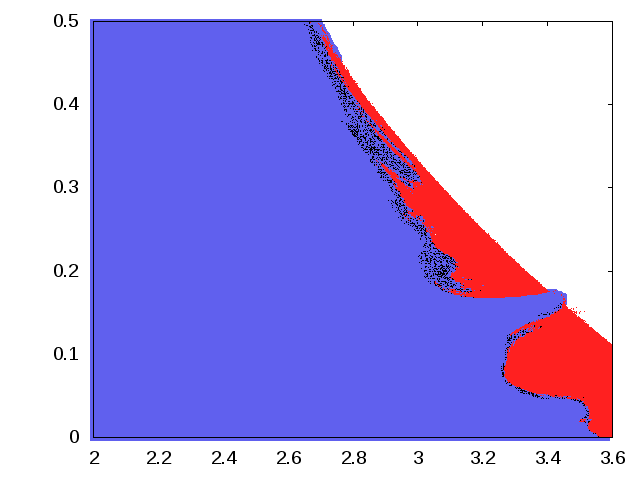}
\includegraphics[width=7.5cm]{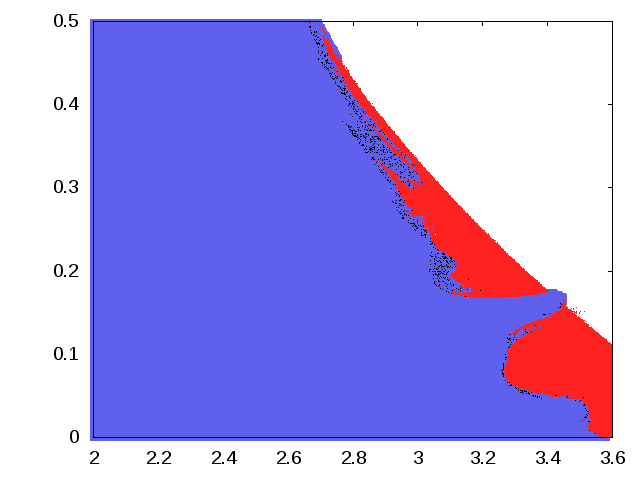}
\includegraphics[width=7.5cm]{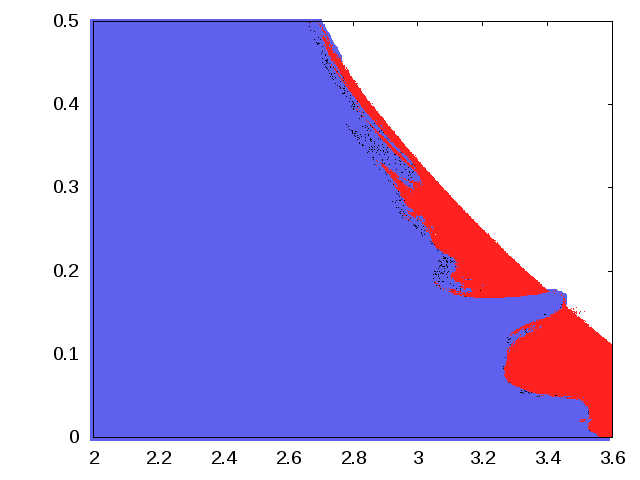}
\caption{Parameter space of the FLM (\ref{FLM}). The axes correspond to
the parameters $\alpha$ (horizontal) and $\eps$ (vertical). Candidates
to SNAs are displayed (in black)  for different orders ($N$) of computation, 
see the text for details. The values of $N$ taken for the computation 
(from left  to right and top to bottom) are $10^4$,
$10^5$, $10^6$ and $10^7$.}
\label{Parameter space diagram 1}
\end{center}
\end{figure}

In this last direction we have noticed that the
detection of SNAs in the parameter space using the
phase sensitive operator (for a description of the method see
\cite{FP95})
%or the appendix \ref{subsection phase sensitivity}
%for a summarized version)
can produce non-reliable results.
In \cite{JT05} it is given an example of a family of
one parametric affine maps in the cylinder
where the supremum norm of the derivative of the invariant curve
goes to infinity when the parameter of the family tends to a certain
critical value, but the invariant curve it still $C^\infty$  before
reaching the critical value. In this case,  the phase sensitive indicator will fail
and it will give a false positive for parameters close to the critical value.
This kind of phenomena can happen when the method is applied to the
FLM giving false positives of the indicator.

To illustrate this statement we have reproduced
the computation of candidates to be SNAs for the FLM done
in \cite{PMR97}. In our computations we have considered several 
values of the parameter $N$ (in the notation of \cite{FP95}) when 
computing the phase sensitivity indicator. We have
fixed a box in the parameter space of the FLM, then  we have discretizated it
in a rectangular grid of points. For each parameter in the grid we have
computed the attractor of the map and its Lyapunov exponent. In the case
of having negative Lyapunov exponent we have computed the phase sensitivity
indicator.

The results are shown in figure \ref{Parameter space diagram 1}. The
parameters  values in white are the ones  for which the iteration of the initial
point diverged to $-\infty$.  The parameters
in red correspond to positive Lyapunov exponent and the ones in
blue correspond to negative exponent.  The points in black correspond to 
the candidates to be SNAs obtained with the phase sensitive indicator. 
The different pictures correspond to
different values of  the order $N$.

We can observe how the number of candidates decays when the
order $N$ of the method is increased. At the end we have much less candidates
 that what appeared to be in the original estimations of \cite{PMR97}.  This is due
to the fact that the indicator can not distinguish between a strange set from
a very wrinkled (but smooth) curve. Then when the number of iterates increases, more and
more candidates to SNAs are discarded.

%%%%%%%%%%%%%%%%%%%%%%%%%%%%%%%%%%%%%%%%%
\subsection{Parameter space and reducibility}
\label{section parameter space and reducibility}
%%%%%%%%%%%%%%%%%%%%%%%%%%%%%%%%%%%%%%%%%5

\begin{figure}[t]
\begin{center}
\includegraphics[width=11cm]{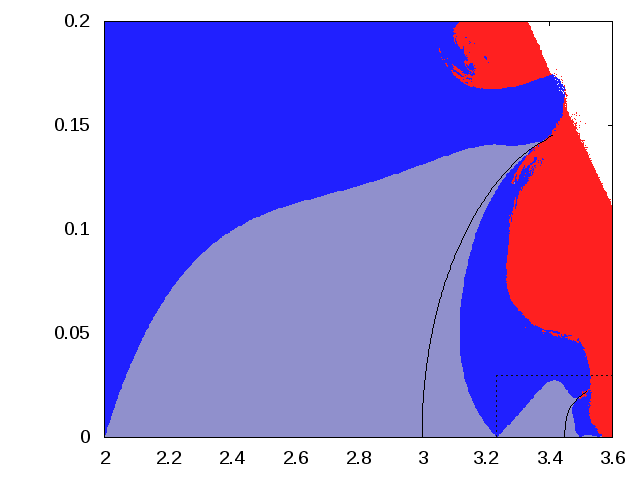}
\includegraphics[width=7.cm]{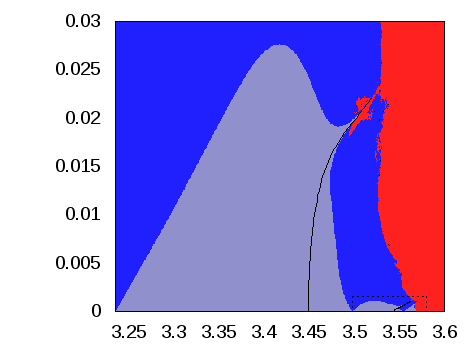}
\includegraphics[width=7.cm]{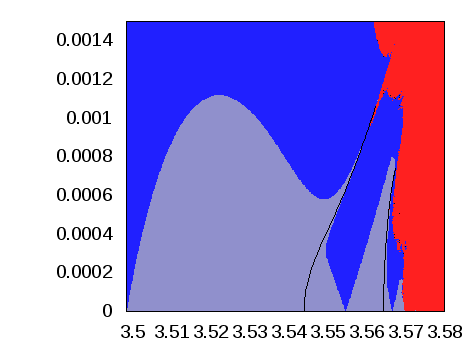}
\caption{Diagram of the parameter space of the FLM (\ref{FLM}) for
$\omega= \frac{\sqrt{5}-1}{2}$. The axis correspond to
the parameters $\alpha$ (horizontal) and $\eps$ (vertical). 
For the correspondence of each color with
the properties of the attractor see table \ref{table coding diagrams}.}
\label{FLM parameter space}
\end{center}
\end{figure}

\begin{table}[t!]
\begin{center}
\begin{tabular}{|lcl|}
\hline
\rule{0pt}{3ex} Color  & \phantom{c} & Dynamics of the attractor \\ \hline
\rule{0pt}{3ex} Black  &\phantom{c} & Invariant curve with zero Lyapunov exponent \\
\rule{0pt}{2ex} Red &\phantom{c} & Chaotic attractor  \\
\rule{0pt}{2ex} Dark blue & \phantom{c} & Non-chaotic non-reducible attractor \\
\rule{0pt}{2ex} Soft blue &\phantom{c} & Non-chaotic reducible attractor \\
\rule{0pt}{2ex} White  &  \phantom{c} & No attractor (divergence to $-\infty$) \\
\hline
\end{tabular}
\end{center}
\vspace{-4mm}
\caption{Color coding for the figure
\ref{FLM parameter space}. }
\label{table coding diagrams}
\end{table}

As in the previous section,  we consider
a certain subset of the parameter space $(\alpha, \eps)$ of the FLM.
The parameters are classified depending on the dynamics of the
 attracting set, but now different properties of the attractor are
considered. There are similar computations 
done in the literature in terms of the Lyapunov exponent of the
attracting set and, in case of having a periodic invariant curve,
its period (\cite{Kan83,PMR98,FKP06}). In our analysis, when we have 
an invariant periodic curve we also take into account whether the 
invariant curve is reducible or not.
Our computation of the zero Lyapunov exponent bifurcation curve has 
been done approximating the invariant curve by its truncated Fourier 
series, instead of approximating the rotation number $\omega$ by rational 
approximations like in \cite{KFP98}. 
 
%The reducibility of an invariant has an important 
%role on the spectral problem associated to the local continuation of the curve, see
%section \ref{section invariant curves in q.p systems}
%or \cite{JT05} for more details.

%As we will see in this section it
%has an important role in the bifurcations 
%diagram of the FLM.

\subsubsection{The procedure}

We have considered three different rectangular subsets in the parameter space and
we have discretized the subset in a grid of points. For each parameter
in the grid we have computed the attracting set by forward iteration. Then
we have computed the Lyapunov exponent of each attracting set. 
To compute this we have approximated the limit (\ref{Lyapunov exponent limit})
for large values of $n$. We stopped when the variation between 
the estimation for $n$, $\frac{3}{4} n$ and $\frac{1}{2}n$ was 
lower than $10^{-3}$.

%(see the appendix \ref{subsection Lyapunov exponents}).
For the parameter values with negative Lyapunov
exponent we assumed that we have a periodic invariant curve and we have
checked if it is reducible as follows. 
Given a initial value $\theta_0\in \T$, we can approximate the value of the
invariant curve $u$  at $\theta_0$ as $u(\theta_0)\approx 
f^M( \theta_0 - M \omega, x_0)$ (for any $x_0$). 
This allows us to compute the attracting curve in a mesh of points. On the other 
hand, we have that (if $f$ is $C^\infty$ and $\omega$ Diophantine) 
an invariant curve is reducible if, and only if, $D_x f(\cdot, u(\cdot))$ has 
no zeros (see Corollary 1 of \cite{JT05}). Then we can use the discrete 
approximation of the curve to check this condition. 

% using the method described in appendix
%\ref{subsection attracting invarian curves}.

The results are shown in figure \ref{FLM parameter space}.
The parameter values where the attracting curve
of the map has zero Lyapunov exponent are also displayed . In practice, this corresponds to the
period doubling of the attracting
invariant curve. The numerical computation of
these bifurcations  has been done as follows.

Given  $u:\T \rightarrow \R$ an invariant curve of (\ref{invariance equation})
we have that it can be approximated numerically by its truncated Fourier series 
at order $N$. For more details on how this can be implemented see 
\cite{CJ00,Jor01}. 
Basically, we compute the Fourier coefficients  as the zero of a suitable 
function $G:\R^{2N+1}\rightarrow\R^{2N+1}$. For the particular case 
of the FLM we can add the parameters $\alpha$ and $\eps$ as unknowns and 
extend $F$ to a function from $\R^{2N+3}$ to $\R^{2N+1}$. 
On the other hand, assume that we have an invariant curve
which is reducible. Then we have that the Lyapunov exponent of
the curve $\Lambda(u)$ given by (\ref{Lyapunov exponent integral})
is indeed differentiable (see \cite{JT05}). If we have a curve $u$ 
with zero Lyapunov exponent, we can use the 
function $\tilde{G} = (G,\Lambda)$ as a continuation 
function for the zero Lyapunov bifurcation curve. Then the bifurcation curve 
in the parameter space is computed as the projection in the $(\alpha,\eps)$ 
coordinates of the set of zeros of $\tilde G$.

To compute this bifurcation curve numerically it can be used a standard continuation
method (see \cite{Sim89}). In our case the kind of continuation that we
do is rather simple. Suppose that the bifurcation curve in the parameter
space is regular (no critical points). Then at every point in the parameter
space the curve can be expressed locally with one of the parameters as a function of
the other. When one of the parameters is fixed, the other can be
computed through a standard Newton method. To do the continuation of 
the bifurcation curve we can vary slightly the fixed parameter and 
compute the corresponding value of the curve. The selection of which 
parameter is fixed has been done depending on the estimated inclination 
of the curve in the parameter space.

\subsubsection{Description of the results}

In figure \ref{FLM parameter space} we show the results of the
computation and we have also plotted
two successive magnifications of certain subsets of the parameter
space. These regions have been marked with a dashed line in the
picture. The points in the parameter space have been codified in
different colors depending on the dynamics of the attractor,
see  table \ref{table coding diagrams}.

At this point let us recall the dynamics of the logistic map. 
We have that for certain values
of $\alpha$ the map has a cascade of period doubling bifurcations.
This is,  we have that there exist a sequence of
values $\{d_n\}_{n\in\N}$,
with $\displaystyle  \lim_{n\rightarrow \infty} d_n =d^*$, such that for
$\alpha \in (d_{k}, d_{k+1})$ the logistic map has an attracting
periodic orbit of period $2^k$ and at the value $\alpha = d_{k+1}$ the
attractor doubles its period (from $2^k$ to $2^{k+1}$). Let us recall
that the logistic map is unimodal for any value of $\alpha$ (i.e. it has
a unique
%turning point which is the 
point where the derivative of
the map is equal to zero). In the period doubling cascade, we also have
 that the attractor crosses the turning point between one period
doubling and the next one. In other words, there exist values
$c_k \in (d_{k}, d_{k+1})$  for which the attracting $2^k$ periodic
orbit of the map is the critical point.

In figure \ref{FLM parameter space} we can observe that from every
parameter value $(\alpha,\eps)= (d_k,0)$ it is born a curve where the
attractor has zero Lyapunov exponent. Let us denote by $D_k$ each of
these curves in the parameter space,  which corresponds to the period 
doubling from period $2^k$ to $2^{k+1}$. In the bifurcation diagram
at the top of figure \ref{FLM parameter space} there are plotted
the curves $D_1$ and $D_2$.
In the bottom left one there are
displayed the curves $D_2$ and $D_3$, and finally in the bottom right
one there are the curves $D_3$ and $D_4$.

In the previous figure it has also been plotted the
reducibility and the non-reducibility regions.
We can also observe that from every parameter
$c_i$ a ``cone of non-reducibility'' is born, in the
sense that there exist two curves in the parameter space which
define a zone where the attracting invariant curve of the map is not reducible.
Let us denote by $C_i^-$ (and respectively by $C_i^+$)
the left (respectively right) boundary of the
non-reducibility region born at the point $c_i$.

We can observe how the curves $C_{i-1}^-$ and $C_{i}^+$ define
an enclosed reducibility region which contains the curve
$D_i$. Moreover these three curves seem to meet in a
tangent way at the same point.

%{\bf A first analysis of the bifurcations diagram} 
%--------------------------------------------------------
\subsubsection{Analysis of the bifurcation diagram}
%--------------------------------------------------------

%fbox{relate the diagram with the truncation of 
%he period doubling cascade.} 

Recall that in section \ref{section basic study of FLM} 
we have illustrated a truncation of the period doubling cascade for a fixed
value of the coupling parameter $\eps$, see figure 
\ref{exemples fractalitzacio a periode 1}. The results described
above on the parameter space of the map agree
with the behavior reported there. Looking at figure
\ref{FLM parameter space} we can observe that
each period doubling is confined inside a reducibility region. Moreover,
when the period is increased each of these regions get closer to the
line $\eps=0$.
Then, if one fixes the value of $\eps$ (arbitrarily small) and let $\alpha$ grow
one should expect a finite number of period doublings. Fixing $\eps$
at a prescribed value $\eps_0>0$ corresponds to fixing a line $\eps=\eps_0$
in the parameter space. If the enclosed regions of reducibility get
closer to $\eps=0$ when the period grows, at some point the
regions will be below the line $\eps=\eps_0$.
On the other hand, the shape of the reducibility regions also explains
why we observe a reducibility loss and afterwards
a reducibility recovery between one period doubling and the next.
At the same time, the bifurcation diagram displays many
interesting phenomena which can be studied.

The first phenomenon that can be observed in the
bifurcation diagram is the birth of non-reducibility cones
around the parameters values $(\alpha,\eps)=(c_i,0)$.
%Note that, using propositions \ref{reduciblity in the analitic case}
%and \ref{normal form near a reduciblility loss}, one can understand the
%loss of reducibility as a bifurcation, in the sense 
%that it is defined as a codimension one condition. 
Recall that (definition \ref{definition of reducibility loss 
bifurcation}) the reducibility loss of an invariant curve
can be seen as a codimension one bifurcation.
This codimension
one condition defines a one dimensional curve in the  parameter space,
which is the boundary
between the reducibility and the non-reducibility regions. With the
notation introduced above, these boundaries will correspond to the
curves $C_i^+$ and $C_i^-$ of the parameter space.

Actually, the existence of a cone of non-reducibility around the
points  $(\alpha,\eps)=(c_i,0)$ is equivalent to
prove that two reducibility loss bifurcations curves are born
at these points. This will be proved under suitable conditions
in \cite{JRT11b},
%chapter \ref{chapter renormalization for q.p. forced maps}
 but the phenomenon can be explained heuristically as follows.

Consider that we have an invariant
(or periodic) curve $u:\T \rightarrow \R$ of the FLM. 
Using corollary 1 of \cite{JT05} we
have that $u$ is reducible if, and only if,
$D_x f_{\alpha,\eps} (\cdot, u(\cdot))$ has
no zeros. For the FLM, we have that the points $(\theta,x)$ for
which $D_x f_{\alpha,\eps}(\theta,x)$ is equal to zero are
those in the set $\{x=1/2\}$. Therefore an invariant curve $x=u(\theta)$ is
reducible if, and only if, it does not intersect the set $x=1/2$.
Respectively, a periodic curve will be reducible if, and only if,
none of its periodic components touches this set.

On the other hand, when  the system is uncoupled
($\eps =0$) we have that the invariant (resp. periodic)
curves are constant lines equal to the fixed
(resp. periodic) points of the uncoupled system. For the
logistic map we have that when the parameter $\alpha$ crosses
the value $c_0$ (resp. $c_i$), the fixed points (resp. one of the
components of the periodic orbit) crosses the set $x=1/2$.
When one adds a small coupling to the system  the invariant
curve is no longer constant. Then when $\alpha$ is increased
crossing the value $c_0$ (resp. $c_i$)
the invariant curve (resp. periodic orbit) has to lose its
reducibility, when it first touches $x=1/2$, and then recover it
again, when it is completely below (or above) $x=1/2$. 
These two contacts
of the invariant (or periodic) curve with
the set $x=1/2$  
give place to the cone of non-reducibility
when the parameters $(\eps,\alpha)$ are close
to $(0,c_0)$ (or $(0,c_i)$). Each boundary of this cone 
correspond to a loss of reducibility
bifurcation curves, namely $C_i^-$ and $C_i^+$, in the parameter space.

%%---------------------------------------

Another observable phenomenon in the bifurcation diagram of the
figure \ref{FLM parameter space} is that from every parameter
value $(\alpha,\eps)=(d_i,0)$ it is born a curve in the parameter
space where the attracting set has
zero Lyapunov exponent. Each of these
curves correspond to a period doubling bifurcation of
the attractor.
%Consider the FLM (\ref{}), using the proposition \ref{prop persistence 
%fixed points} we have that any hyperbolic period orbit of the 
%logistic maps for $\alpha\in (f_k, f_{k+1})$ (with  $k>0$)  
%will become periodic invariant curves of the FLM for $\eps$ 
%sufficiently small. On the other hand, 
For diffeomorphisms it is known (\cite{BHTB90}) that the period
doubling of an invariant torus is a codimension one bifurcation. Recall
that the FLM is not invertible, therefore this theory is
not directly applicable.
Nevertheless, in the reducible case, one can apply a normal form
procedure around the invariant curve, to obtain a one
dimensional system in a neighborhood of the curve. We will formalize this
argument in section \ref{section period doubling and reducibility}.
Then we can understand the period doubling
bifurcation of the invariant curve as a period doubling of the
reduced system (in the classical one dimensional sense). For the
uncoupled case ($\eps=0$) this happens at each parameter value
$\alpha=d_i$. Due to the codimension one character of the bifurcation
a period doubling  bifurcation curve is born at each of
the parameters $(\alpha,\eps)= (d_i,0)$.

% -------------------------------- 

% -------------------------------- 

\begin{figure}[t]
\begin{center}
\includegraphics[width=7.5cm]{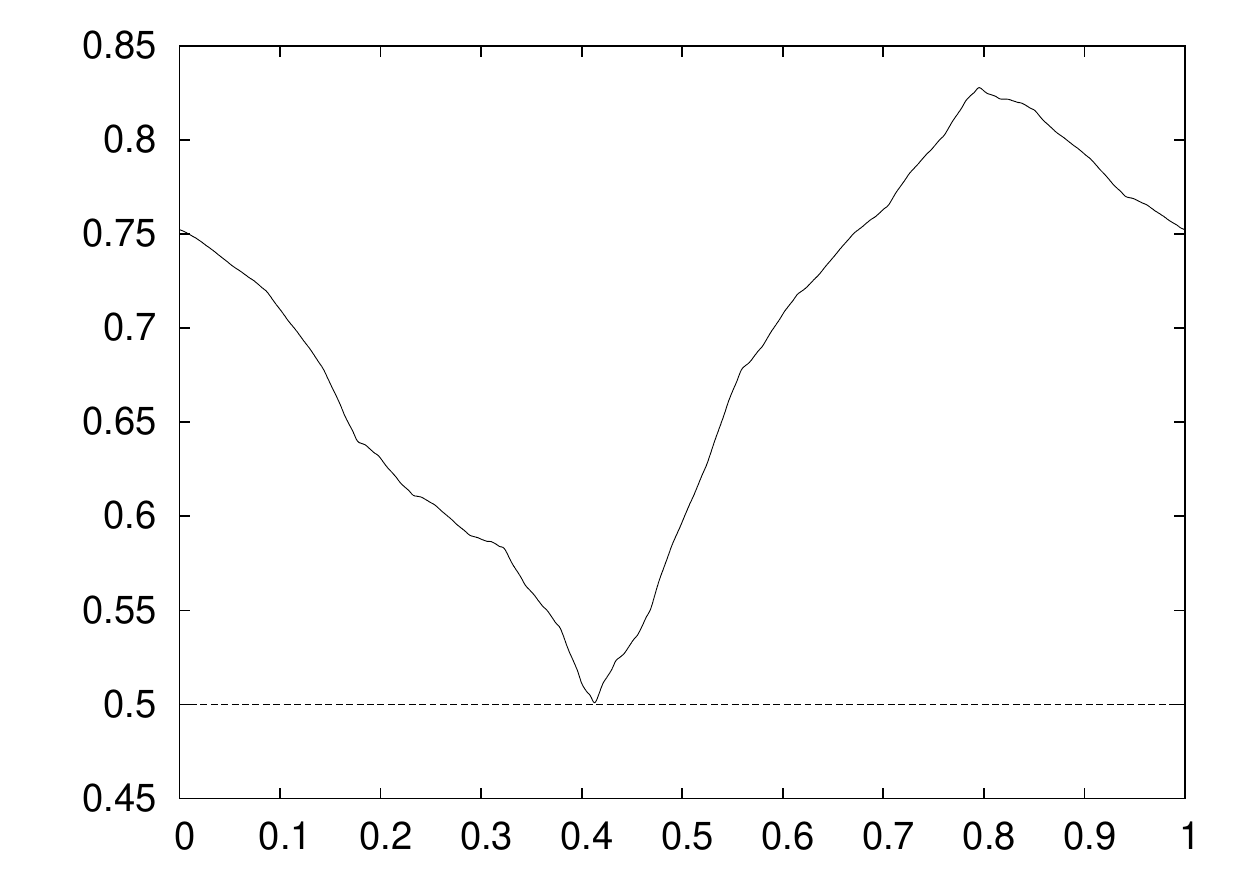}
\includegraphics[width=7.5cm]{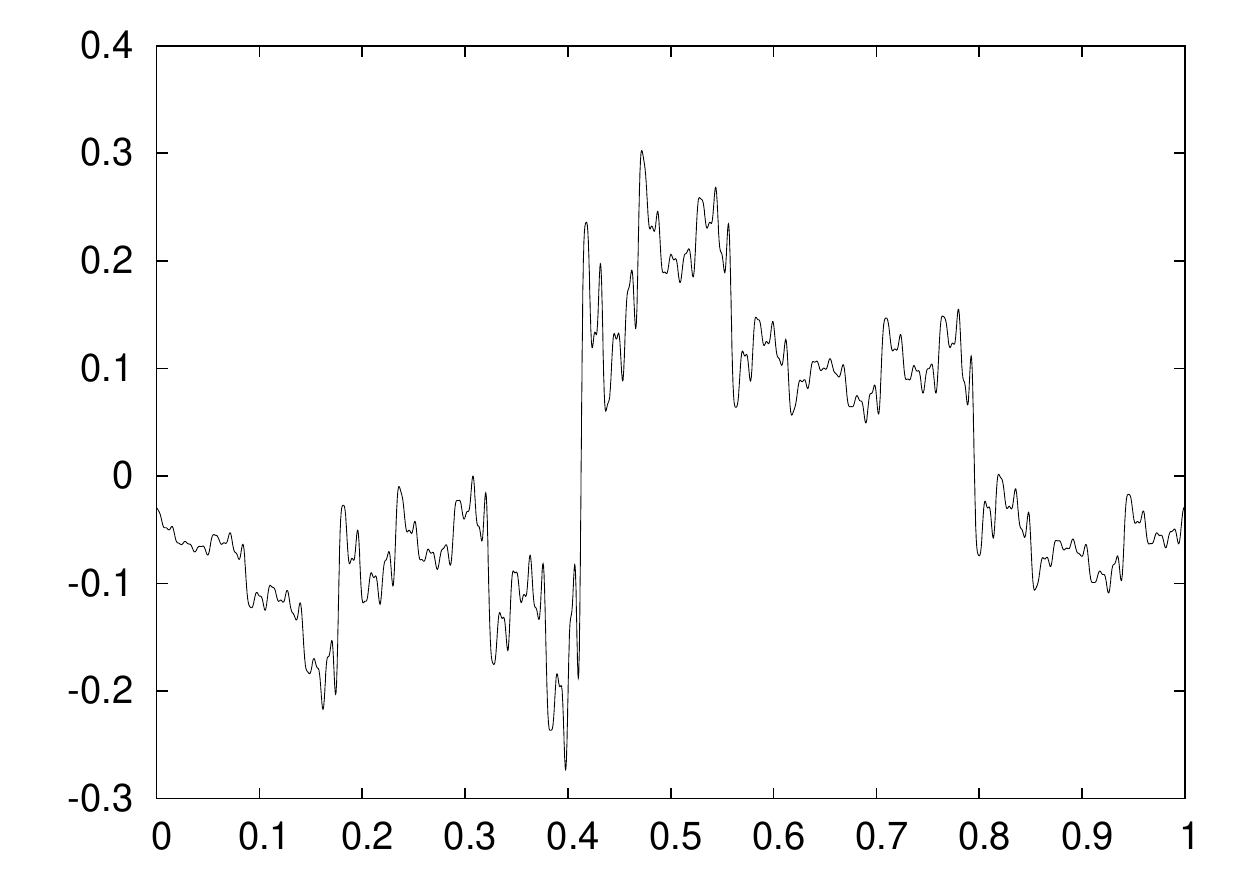}
\includegraphics[width=7.5cm]{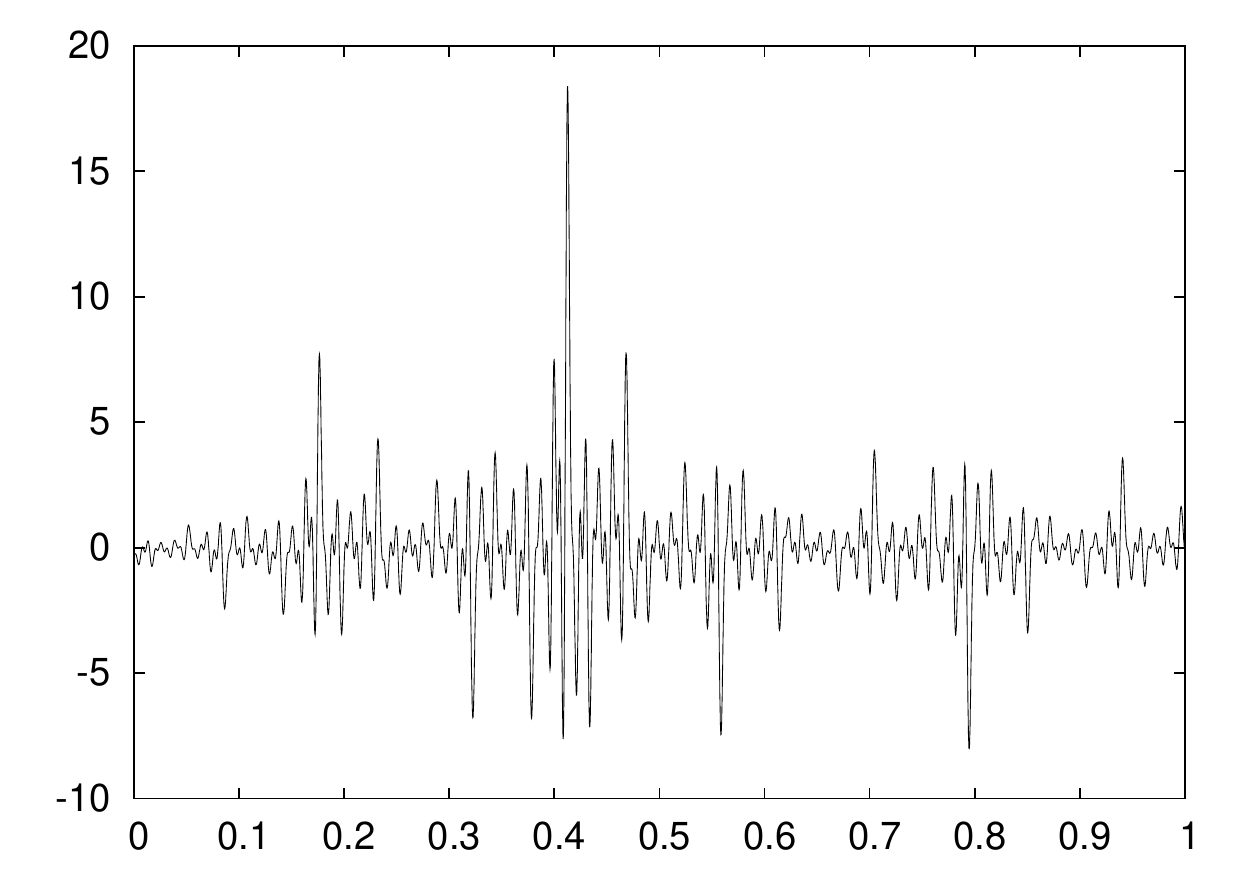}
\includegraphics[width=7.5cm]{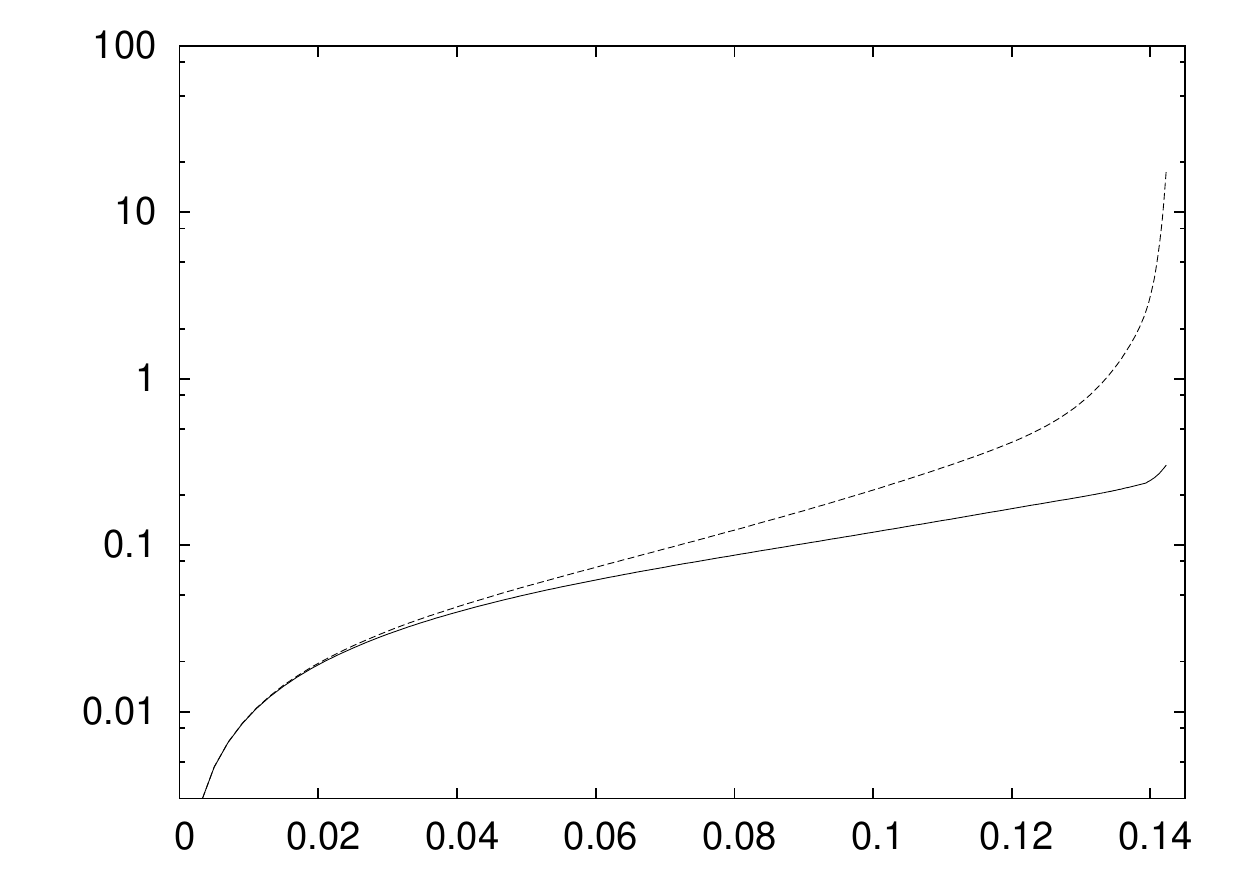}
\caption{In the top left we have the invariant curve $x(\theta)$ with
zero Lyapunov exponent of the FLM (\ref{FLM}) for the terminal point of 
the first period doubling bifurcation. In the top right and bottom left 
we have the first ($x'(\theta)$) and the second derivative ($x''(\theta)$)  
of the same invariant curve. In these pictures the 
horizontal axis corresponds to $\theta$. 
In the bottom right we have (in the vertical
axis on a logarithmic scale)  
$\displaystyle \sup_{\theta\in\T} |x_\eps'(\theta)|$ (solid
line) and $\displaystyle \sup_{\theta\in\T} |x_\eps''(\theta)|$ (dashed line) as a graph
of the parameter $\eps$ (horizontal axis),
where $x_\eps$ is the invariant curve with zero Lyapunov exponent for each
$\eps$.}
\label{critical curve}
\end{center}
\end{figure}

%%--------------------------------------

%%--------------------------------------
One remarkable property (observed numerically)  is that each period doubling 
bifurcation curve $D_i$ is confined in a reducibility region delimited 
by $S_{i-1}^+$ and $S_{i}^-$.
This can be explained analyzing the numerical procedure.
Recall that the continuation of the bifurcation curve has
been done as long as the estimated error was below a
prescribed tolerance. On the other hand the
function used for the numerical continuation  contains the Lyapunov
exponent of the invariant curve as one of its
components. % (see appendix \ref{section computation bifurcation curves}).
The method used to estimate the Lyapunov exponent behaves much
worse when the reducibility is lost. Note that the Lyapunov 
exponent is obtained integrating numerically 
$\ln|D_x f(\theta, u(\theta))|$. In the $C^\infty$ case
we have have that, in the reducible case the function is $C^\infty$ 
whereas in the non-reducible case the function is not even bounded. 
This makes the numerical integration behave much worse when the reducibility 
is lost. 
% (see appendix
%\ref{subsection Lyapunov exponents}). 
%This change
%of behavior makes the method fail when the reducibility is lost.
Despite of that, we have the following reasons to believe that the invariant
curve (at the parameter of bifurcation) is actually destroyed due
to the reducibility loss. 

%\begin{enumerate}
%\item
From the theoretical point of view, when the reducibility is
lost the spectral problem associated to the continuation
of the invariant curve changes drastically. When one tries to apply 
the IFT to an invariant curve there is a big difference from the 
reducible case to the non-reducible one (\cite{JT05}). 
%(see the comments on the proof of theorem \ref{continuation theorem}
%in section \ref{subsection continuation of invariant curves}).
Given an invariant curve and its skew product consider the transfer
operator %(\ref{transfer operator}) 
associated to the curve. We
have that the invariant curve can
be continued when $1$ does not belongs to the spectrum of the transfer
operator. In the reducible case the spectrum corresponds to the
circle of radius equal to the exponential of the Lyapunov exponent of
the curve. In the non-reducible case the spectrum corresponds to the 
disk of the same radius. When a period doubling occurs, there is an 
attracting invariant curve which becomes unstable, in other 
words its Lyapunov exponent crosses zero. In the reducible case, 
one might expect the spectrum  of the transfer operator 
associated to the invariant curve (which is a circle) 
to cross $1$ transversally when the bifurcation happens, the curve 
becomes unstable but survives. 
But in the non-reducible case there is no chance of transversal cross. 

%  ------------------------------------------- 
%\item 
From the numerical point of view, one can study the
behavior of the invariant curve for the parameter values
in the period doubling bifurcation just
before losing the reducibility. In figure \ref{critical curve}
we have displayed the invariant curve $x(\theta)$ of the
FLM (\ref{FLM}) for the parameters $(\alpha,\eps) = (3.3796, 
0.1423)$. This parameter corresponds to the last point in the
invariant curve that we have obtained in the numerical computation.
In the figure we also show the first and second derivatives
of the curve.

We can observe that the invariant curve $x(\theta)$ seems to keep
continuous and it does not fractalize\footnote{
Heuristically, we consider that a curve fractalizes if its length
goes to infinity.}  but what fractalizes is its
derivative $x'(\theta)$. Back to figure \ref{FLM parameter space}
we can observe that the period doubling
bifurcation curve can be parameterized by the parameter $\eps$.
Then for each $\eps$ we can consider $\alpha(\eps)$
the parameter for the period doubling bifurcation curve
 and  $x_\eps(\theta)$ the invariant curve with zero Lyapunov
exponent for the corresponding values of $\alpha$ and $\eps$.
In figure \ref{critical curve} we have the supremum
of $|x'_\eps(\theta)|$ and $|x''_\eps(\theta)|$ as a
graph of $\eps$. The figure indicates that $\displaystyle \sup_{\theta\in\T} 
|x''_\eps(\theta)|$ grows unbounded when we approximate certain
critical value $\eps \approx 0.1423$. Moreover, looking at
the graph of $|x''_\eps(\theta)|$ this limit seems to be uniform
when taking the supremum on any subinterval of $\T$.

This destruction
process is similar to the fractalization process described
 in \cite{JT05}, but the curve which fractalizes here is
$x'(\theta)$. The fact that the end point of the period
doubling bifurcation corresponds to the parameter where
the invariant curve touches the line $x=1/2$ was already observed
 in \cite{KFP98}. The computations of the bifurcation curve
done in the cited work are based on the rational approximation of the
rotation number $\omega$. Our method
is based on the computation of the Fourier series of the
invariant curve and it has the advantage  that it allow us
to compute the derivatives of the curves easily.

%\end{enumerate}

To finish this section let us remark that this analysis is 
far from complete. In this direction, in the forthcoming sections
\ref{section obstruction to reducibility} and
\ref{section period doubling and reducibility}
we  present two studies in order to 
understand better the  bifurcation diagram of 
the figure \ref{FLM parameter space}. 
In section \ref{section summary and conclusions} we 
summarize the different results and we analyze their implications 
on the cited bifurcation diagram.

%%%%%%%%%%%%%%%%%%%%%%%%%%%%%%%%%%%%%55
\section{Obstruction to reducibility} 
\label{section obstruction to reducibility}
%%%%%%%%%%%%%%%%%%%%%%%%%%%%%%%%%%%%%55

Recall that the set $x=0$ is always an invariant curve of the FLM. 
Apart from this set, when $\eps=0$ and $\alpha>1$ we have also the invariant 
curve given as $x_{\alpha,0}(\theta)= 1 - \frac{1}{\alpha}$. Using proposition 
\ref{prop persistence fixed points} we have that for any $\alpha \neq 3$ 
this curve persists for small values of $\eps$. Let us denote by 
$x_{\alpha,\eps}$ the continuation of this curve. In this section we are 
concerned on the reducibility of this invariant curve. %More concretely 
Considering the images and preimages of the set of points 
where the derivative of the maps is equal to zero, we 
will construct regions in the parameter space where 
the curve $x_{\alpha,\eps}$ can not be reducible. For instance 
we will see that  $\alpha>2$ then  $|\eps| < 1- \frac{1}{2}$ is 
a necessary condition for the reducibility of $x_{\alpha,\eps}$.

In this section we will consider different sets in the cylinder 
$\T \times [0,1]$. These sets will be the closed graph of 
a curve or a subset of a graph. In general when we say that one of these
sets is above (respectively below) of another, we 
mean that, for each value of $\theta \in \T$ the corresponding $x$-coordinate 
of the first set is bigger (resp. smaller) than the $x$-coordinate of the 
other set (for the same $\theta$).  The proofs in this section 
have been omitted because of  their  simplicity. 

%%%%%%%%%%%%%%%%%%%%%%%%%%%%%%%%%%%%%%%%
\subsection{A first constraint on the reducibility}
%%%%%%%%%%%%%%%%%%%%%%%%%%%%%%%%%%%%%%%%

\begin{defin} 
Given a q.p. forced map like (\ref{q.p. forced map general}), we define 
its  {\bf critical set} as the set of points on its domain 
where the derivative of the map (with respect to  $x$) is zero. In other words, 
\[
P_0 = \{
(\theta, x) \in \T\times \R | \thinspace D_x f(\theta,x) =0
\}.
\]
\end{defin}

In the case of the logistic map we have that $P_0=\{(\theta, x) 
\in \T\times \R | \thinspace x=1/2\}$. 

When we have a q.p. forced map like (\ref{q.p. forced map general}) which
 is $C^\infty$ and $\omega$ is Diophantine, Corollary 1 of \cite{JT05} implies 
that an invariant 
curve $y_0(\theta)$ is reducible if, and only if,
\[
P_{0} \cap \Graph (y_0)=\emptyset,
\]
where $\Graph (y_0) =  \{(\theta, x) \in \T\times [0,1] |\thinspace  
x=y_0(\theta) \}$.

In the particular case of the FLM when 
$\alpha(1+|\eps|)<4$ we have that the compact cylinder $\T\times [0,1]$ is 
invariant by the map. Then, any reducible invariant curve is 
either above or below the critical set. Actually, when the map 
is uncoupled ($\eps=0$) we have that the invariant curve $x_0(\theta) = 1-1/\alpha$ is 
above the critical set when $\alpha > 2$.

Consider now $P_1$ the image of the critical set by the map, namely 
{\bf post-critical set}, which is defined as
\[
P_1 = \{(\bar{\theta}, \bar{x}) \in \T\times [0,1] | 
\thinspace \bar{\theta} = \theta +\omega, \text{ }
\bar{x}= f(\theta,x), \text{ for some } (\theta,x)\in P_0\}.
\]
In the particular case of the FLM we have, 
\[
P_1 = \left\{(\bar{\theta}, \bar{x}) \in \T\times [0,1] | 
\thinspace \bar{x}= \frac{\alpha}{4}(1+\eps \cos(2\pi(\bar{\theta} - \omega)))\right\}.
\]

\begin{prop}
\label{proposition postcritical set} 
In the case of the FLM we have the following properties on the 
post-critical set when $|\eps|<1$.
%\begin{proposition}
\begin{enumerate}
\item The set $P_1$ is above the image of any other point 
$(\theta,x) \in \T\times [0,1]$. 
\item The points below $P_1$ have two preimages in $\T\times \R$, the points 
in $P_1$ have one preimage and the points above $P_1$ have no preimage. 
\end{enumerate}
\end{prop}

Consider the invariant curve $x_{\alpha,\eps}$ of the FLM 
introduced at the beginning of this section. When $\eps=0$ and 
$\alpha>2$ we have that the invariant curve exists and it is above 
the critical set $P_0$. On the other hand, $x_{\alpha,\eps}$ is always
below $P_1$ due to proposition \ref{proposition postcritical set}. 
Concretely, we have that a necessary condition for the reducibility 
of $x_{\alpha,\eps}$ is that $P_0$ and $P_1$ do not intersect. 
Otherwise there is no room for the invariant curve to exist without 
losing its reducibility. 

Then for any $\alpha>2$ it is necessary that,
\[
\frac{1}{2} < \frac{\alpha}{4}( 1 +\eps \cos(2\pi\theta)) \text{ for any } \theta \in \T,
\]
which will be satisfied if, and only if, 
\begin{equation}
\label{A first bound for the reducibility}
|\eps| < 1 - \frac{2}{\alpha}. 
\end{equation}

This gives a first constraint in the parameter space for the reducibility of
the curve $x_{\alpha,\eps}$. In figure \ref{bounds of reducibility} 
we have displayed this constraint together with the bifurcation diagram 
computed in section \ref{section parameter space and reducibility}.
We can see that the constraint is quite sharp for values of $\alpha$ close to 
$2$, but it becomes rapidly pessimistic for values of $\alpha$  far from $2$. 

%%%%%%%%%%%%%%%%%%%%%%%%%%%%%%%%%%%%%%%%%
\subsection{Further constraints on the region of reducibility}
%%%%%%%%%%%%%%%%%%%%%%%%%%%%%%%%%%%%%%%%%

\begin{figure}[t]
\begin{center}
\includegraphics[width=7.5cm]{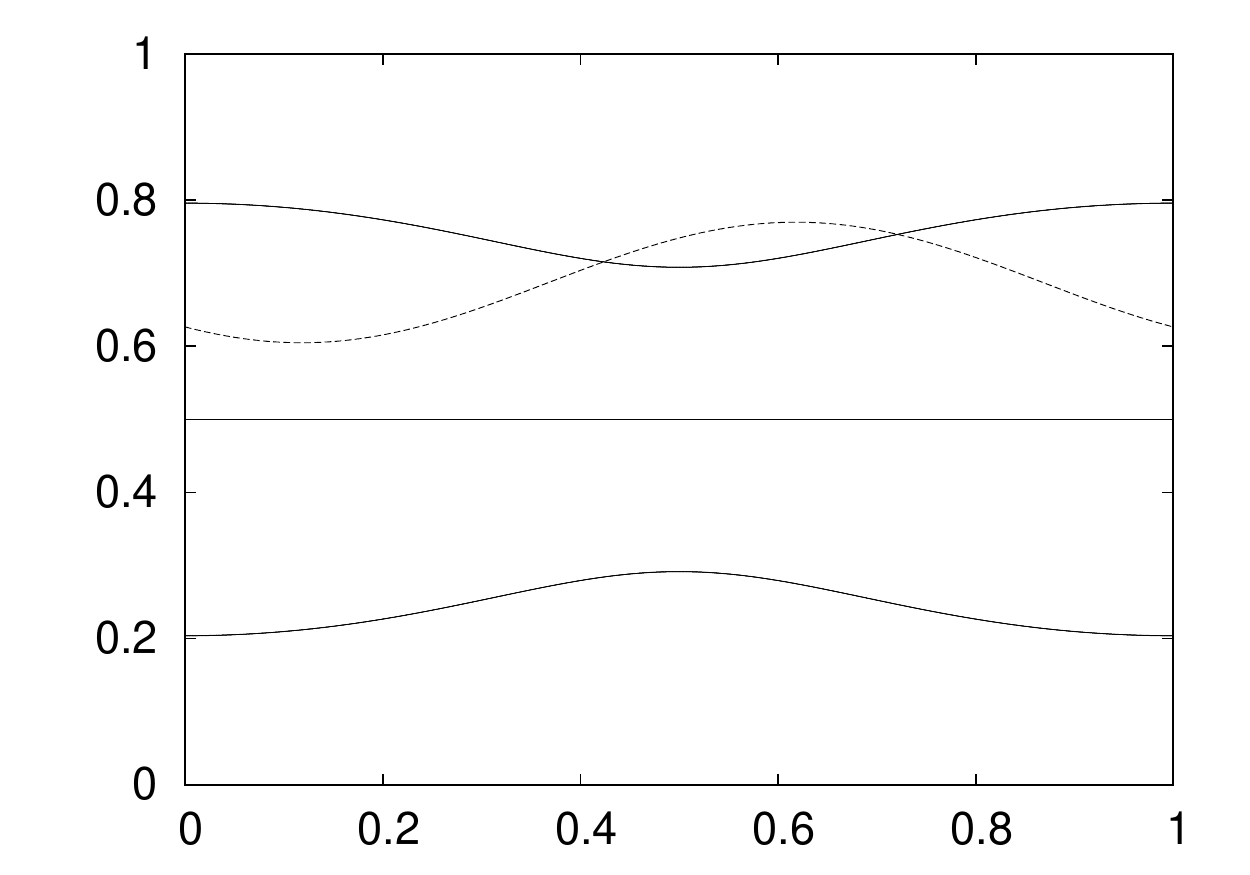}
\includegraphics[width=7.5cm]{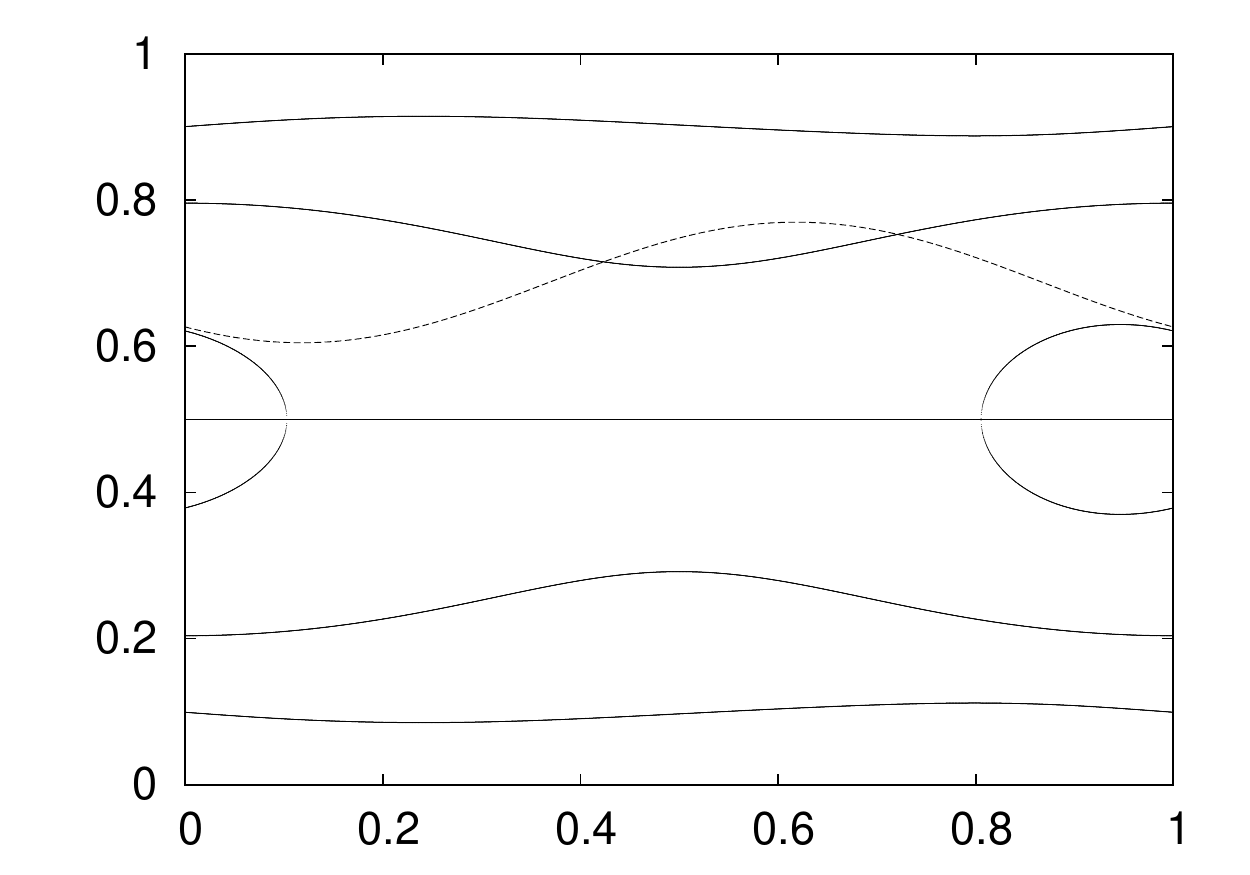}
\includegraphics[width=7.5cm]{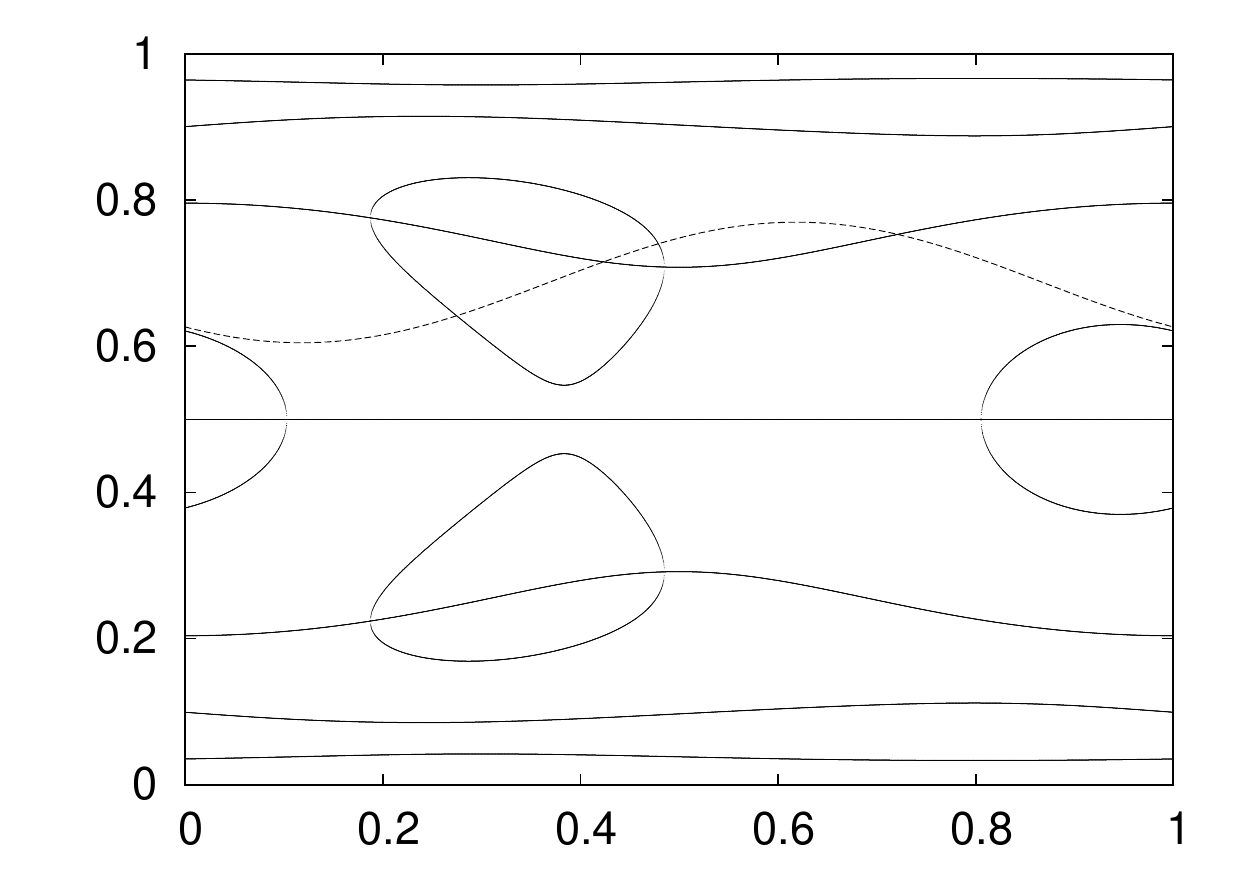}
\includegraphics[width=7.5cm]{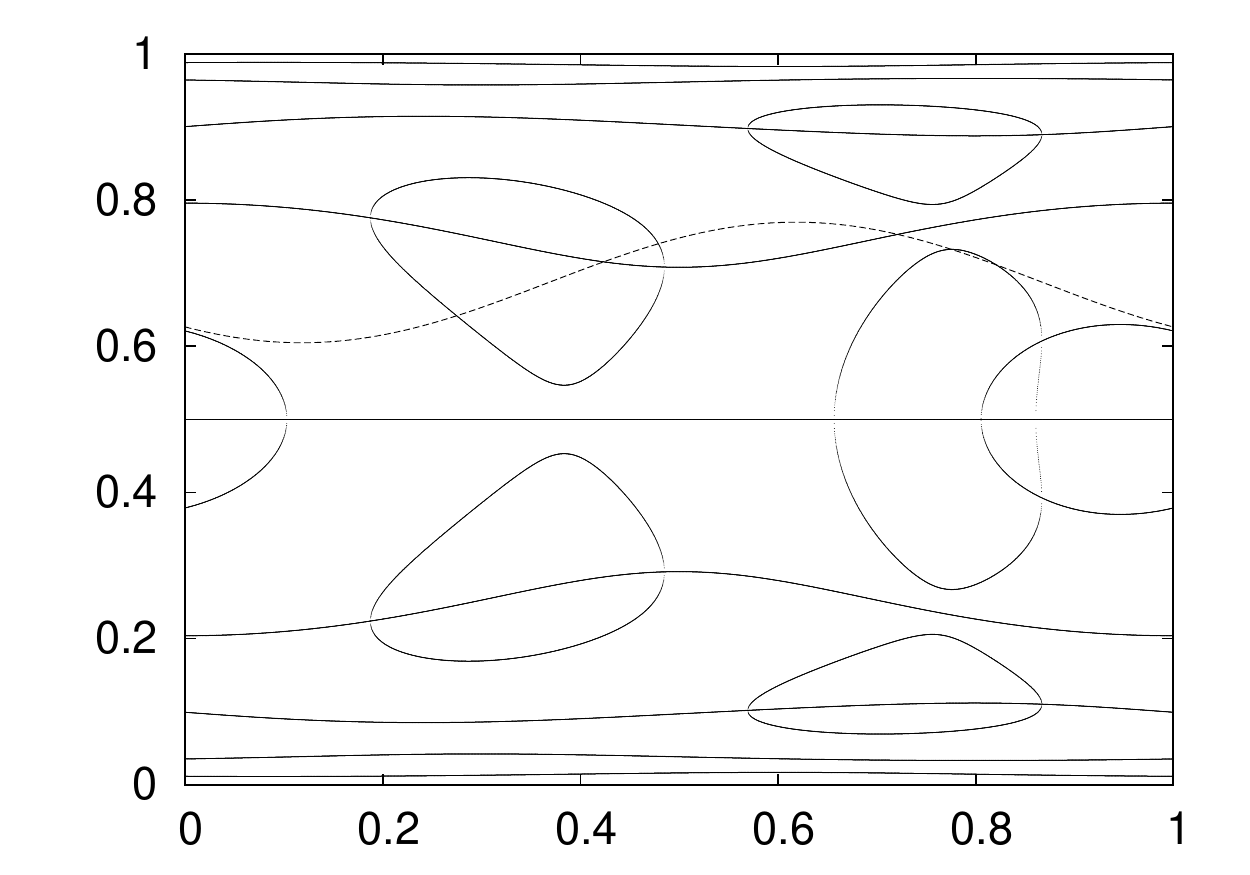}
\caption{ The union of the pre-critical set of the FLM for 
$(\alpha,\eps)=(2.75,0.12)$, from top to bottom and left to right 
we have the union of the $k$ first pre-critical sets, for $k=1,2,3,4$.
The horizontal axis corresponds to $\theta$ and the vertical one to $x$.} 
\label{preimages of the critical line}
\end{center}
\end{figure}

In this section, we consider the preimages of the critical 
set to give  additional restrictions to the region of the 
parameter space where the invariant curve 
$x_{\alpha,\eps}$ can be reducible. 

\begin{defin}
Given $F$ a q.p. forced map like (\ref{q.p. forced map general}) we define 
$P_{-k}$ the $k$-th {\bf pre-critical set} as the $k$-th preimage of the 
critical set $P_0$. In other words,  
\[
P_{-k} = \{(\theta, x) \in \T\times \R |\thinspace  F^k(\theta,x) \in C_0 \}. 
\]
\end{defin}

In figure \ref{preimages of the critical line} we have (in a solid line)
the first four pre-critical sets for the FLM for a concrete values of
the parameters. In the figure we have (in a dashed line) 
the post-critical set 
$P_1$. In the top left part of the figure we have displayed the sets 
$P_0$ and $P_{-1}$. Note that $P_0$ lays completely below the set $P_1$, 
therefore using proposition \ref{proposition postcritical set} we have that 
each point of $P_0$ has two preimages, which will belong to $\T\times [0,1]$. 
These two preimages
 form the set $P_{-1}$. In the top right part of the figure we 
have added the set $P_{-2}$ to the previous ones. The lower component of 
$P_{-1}$ is below $P_1$ (so it has two preimages), which form part of 
the set $P_{-2}$. The upper component of $P_{-1}$ intersects $P_1$. Then only 
some part of it has preimage in $\T\times [0,1]$. Indeed, the preimages form two
symmetric arches around $P_{0}$. Note that the preimage of $P_1$ is $P_0$, 
then the preimages of the arch of $P_{-1}$ connecting 
two points of $P_0$ are two arches connecting two points of  $P_0$. 

When further preimages are considered,  higher pre-critical sets are obtained. 
The number of components and the shape of these components depend on the relative 
position of the previous pre-critical set. For example in the bottom part of the 
figure \ref{preimages of the critical line} we have added the sets 
$P_{-3}$ (left) and $P_{-3} \cup P_{-4}$ (right). 
We can observe that the two arches of $P_{-2}$ around 
$P_0$ give place to two pairs of arches of $P_{-3}$, one 
around each component of $P_{-1}$. 
Some parts of the two arches of $P_{-3}$ lay below  $P_1$, then when we consider 
their preimage they give place to new arches of $P_{-4}$ 
around $P_{0}\cup P_{-2}$.

We are considering these pre-critical sets because they suppose
an obstruction to the reducibility of the invariant curve $x_{\alpha,\eps}$. Assume 
that  $y_0$ is an invariant curve by the map. Since its 
graph is invariant by the map we have 
that if it intersects $P_{-k}$ for some $k$ it 
will eventually intersect the set $P_0$. Therefore given an invariant 
curve $y_0$, a necessary condition for its reducibility is that 
\[
P_{-k} \cap \Graph (y_0)=\emptyset \quad \forall k \geq 0,
\] 
where $\Graph (y_0) =  \{(\theta, x) \in \T\times \R |\thinspace  
x=y_0(\theta) \}$.

For example, the FLM for the values $\alpha$ and $\eps$ of the 
figure \ref{preimages of the critical line} we can see that 
there is a component of $P_{-4}$ which connects $P_0$ with $P_1$. 
Then we have that it can not exists an invariant reducible 
curve above $P_0$.

We want to formalize this discussion on the obstruction to the reducibility to 
give additional restrictions on the reducibility of the set $x_{\alpha,\eps}$ 
the invariant set of the FLM  introduced at the beginning of this
section. To do that we need to 
analyze with some more detail the pre-critical set.
We have the following improvement of the proposition 
\ref{proposition postcritical set}.

\begin{prop}
\label{proposition preimages FLM} 
Consider the FLM  with $\alpha(1+|\eps|)<4$ and the post-critical set 
$P_1 =  \left\{(\bar{\theta}, \bar{x}) \in \T\times \R | 
\thinspace \bar{x}= \frac{\alpha}{4}(1+\eps \cos(2\pi(\bar{\theta} - \omega)))\right\}$. 
Consider also $S_1$ the set of points of the compact cylinder 
 $\T\times [0,1]$ which are below or in
$P_1$, in other words 
\[
S_1= \left\{(\theta,x) \in \T\times [0,1] \left| \thinspace 
x \leq \frac{\alpha}{4}(1+\eps \cos(2\pi(\bar{\theta} - \omega))) \right.  \right\}.
\]
Then we have that any point $(\theta,x)\in S_1\setminus P_1 $  
has exactly two preimages which are given by $H_+(\theta, x)$ and 
$H_-(\theta,x)$, where 
\[
%\begin{equation}
%\label{}
\begin{array}{rccc}
H_\pm: & S_1 & \rightarrow & \T^1 \times [0,1]  \\
& \left( \begin{array}{c} \theta \\ x \end{array} \right) &
\mapsto &  \left( \begin{array}{c}  \theta - \omega \\ 
\rule{0ex}{6ex}\displaystyle \frac{1}{2} \mp \sqrt{\frac{1}{4} - 
\frac{x}{\alpha(1 +\eps \cos(2\pi(\theta - \omega)))}}
\end{array} \right)
\end{array}
%\end{equation}
\]
Moreover we have that $H_+$ (respectively $H_-$) maps homeomorphically the set 
$S_1$ to the set $\T\times [0,1/2]$ (and respectively $\T \times[1/2,1]$). 
Finally we have that the map $H_+$ preserves the orientation, in the sense that, it
is monotone with respect to $x$.  
On the flip side the map $H_-$ reverses orientation, i. e. it 
swaps relative positions (with respect to the $x$-coordinate). 
\end{prop}

\begin{figure}[t]
\begin{center}
\includegraphics[width=12cm]{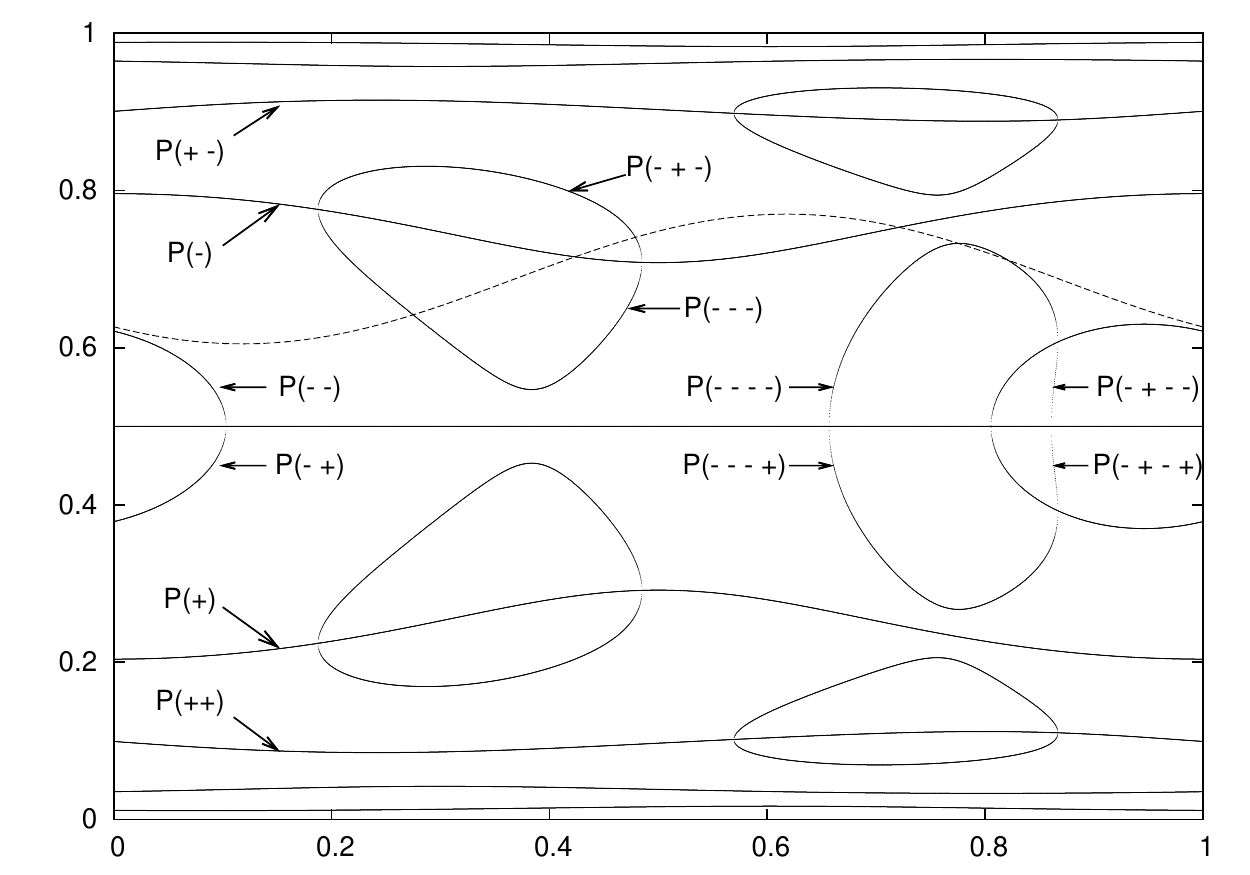}
\caption{ An example for the symbolic codification of the pre-critical 
sets. We have the union of the pre-critical set $P_0$ to $P_4$, and 
the post-critical set $P_1$ of the $FLM$ for 
$(\alpha,\eps)=(2.75,0.12)$. We have indicated the symbolic codes of 
some of the components of the pre-critical sets. The horizontal axis 
corresponds to $\theta$ and the vertical one to $x$.}
\label{symbols of precritical lines}
\end{center}
\end{figure}

%Going back to the original problem, we want to gives bounds in the parameters 
%space for the parameter values where the curve $x_{\alpha,\eps}$ is reducible. 
%We have already seen that when $\alpha>2$ and $\eps=0$ the curve 
%is above the set $P_0$, and then we have seen that the bound 
%(\ref{A first bound for the reducibility})
%|\eps| < 1 - \frac{2}{\alpha}. 
%is a necessary condition to have  $x_{\alpha,\eps}$ a reducible curve. 

This result can be used to describe the pre-critical set with 
some more detail. By construction we have that when  $\alpha>2$ 
the constraint  
(\ref{A first bound for the reducibility}) is satisfied, consequently the set $P_{1}$ is 
strictly above $P_0$. We can apply the proposition 
\ref{proposition preimages FLM} to obtain that 
the set $P_{-1}$ is composed by the union of two different sets,  
i. e. $ P_{-1}= P(+) \cup P(-)$, where 
\[
P(+) = H_+(C_0) \text{ and } P(-)=H_-(C_0).
\]
Moreover we have that $P(+)$ is below $P_0$ and $P(-)$ is above it. 

Now we can consider further preimages. Since $P(+)$ is below $P_0$ we have 
that it belongs to $S_1$ therefore its preimages are defined. Let us denote 
by $P(++) = H_+(P(+)) = H_+ \circ H_+ (P_0)$ and $P(+-)= H_-(P(+)) 
= H_-\circ H_+ (P_0)$. 

On the other hand when,  we consider the preimages of $P(-)$, 
we can have different relative positions 
between $P(-)$ and $P_1$ depending on the parameters. It might happen  
that the curve $P(-)$ is completely below $P_1$, completely above it
or that they intersect. In the case that $P(-)$ has points above 
$P_1$ we have that $H_\pm$ is not defined in these points, therefore 
the set $H_\pm(P(-))$ is not well defined. This can be fixed 
if we formally extend $H_\pm$ to $\T \times [0,1]$ as 
\[
H_\pm (\theta,x) = \left\{ 
\begin{array}{ll} 
H_\pm (\theta,x) & \text{ if } (\theta,x) \in S_1 \\
\emptyset & \text{ otherwise}. 
\end{array}
\right.
\]
With this extension we  can consider the sets $P(-\pm) = H_\pm(P(-))$ without 
problems of definition. If the set $P(-)$ is completely above $P_1$ 
we will have that $P(-\pm) = \emptyset$, and if the set $P(-)$ is partially 
above $P_1$ we have that $P(-\pm)$ is the 
preimage of the points below $P_1$. With this notation we have that 
\[
P_{-2}= P(++) \cup P(+-) \cup P(-+) \cup P(--).
\]

This symbolic codification of the different components of $P_{-2}$ has a 
straight forward generalization to the sets $P_{-k}$. 

\begin{prop} 
Let $\{+,-\}^k$ denote the Cartesian 
product $k$ times of the set of two elements $\{+,-\}$. 
Given $s=(s_1, \dots, s_k) 
\in \{+,-\}^k $, let us define the set $P(s)$ as 
\[
P(s) = H_{s_k} \circ H_{s_{k-1}} \circ \dots \circ H_{s_1} (P_0). 
\]

Then we have,
\[
\displaystyle P_{-k} = \bigcup_{s\in \{ -,+ \}^k} P(s). 
\]
\end{prop}

In figure \ref{symbols of precritical lines} we have the
four first pre-critical sets of the FLM for the parameter values
$(\alpha,\eps)=(2.75, 0.12)$. 
The symbolic codification of some of the components of the 
pre-critical set have been indicated in the picture. The different 
sequences can be deduced using proposition 
\ref{proposition preimages FLM}.

Now we will see which components of the sets $P_{-k}$ can pose 
an obstruction to the reducibility of the invariant curve $x_{\alpha,\eps}$.

\begin{figure}[t]
\begin{center}
\includegraphics[width=12cm]{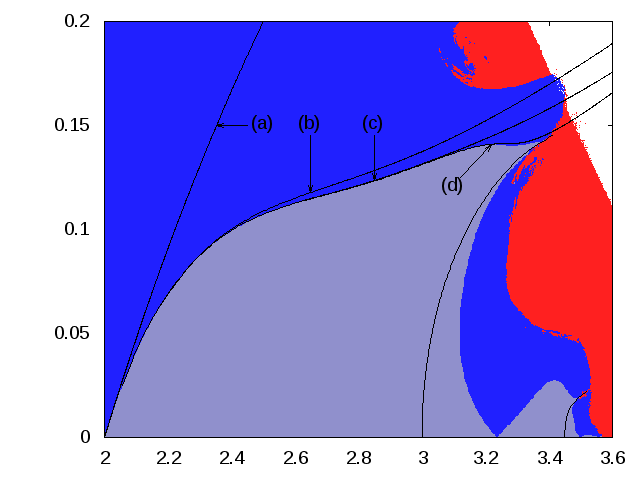}
\caption{The different constrains on the region of 
reducibility of the invariant curve 
$x_{\alpha,\eps}$ together with the bifurcation diagram of section 
\ref{section parameter space and reducibility}. 
The curve (a) corresponds to the constraint 
(\ref{A first bound for the reducibility}), the curves (b), (c), (d) are the 
constrains  associated to the 
pre-critical set $P((-)^2)$, $P((-)^4)$ and $P((-)^{12})$ respectively.}
\label{bounds of reducibility}
\end{center}
\end{figure}

Consider the pre-critical sets of the FLM for $\alpha>2$ and suppose that the 
first constraint of reducibility (\ref{A first bound for the reducibility})
is satisfied. Then we have that $P(-)$ 
is a closed curve above $P_0$. To construct 
additional constrains on the reducibility let us assume that $P(-)$ and $P_1$ intersect 
in exactly two different points. If we consider the preimages of 
$P(-)$ we have that $P(-+)$ is connected by an arch to $P_0$ and below it. On 
the other hand $P(--)$ is also an arch connected to $P_0$ but above it. Now 
it can happen that $P(--)$ intersects the post-critical set $P_1$. If such is 
the case then we would have that one of the component of $P_{-2}$ connects
 $P_0$ with $P_1$. This would imply that no invariant reducible curve can exists
 between $P_0$ and $P_1$.

To avoid the situation described above, one might require the set 
$P(--)$ (when defined) to be below the set $P_1$. This condition gives place 
to an additional constraint in the parameter space of the FLM. In figure 
\ref{bounds of reducibility} we have this constraint on the reducibility 
together with the first constraint and the bifurcation diagram of section 
\ref{section parameter space and reducibility}.

When the set $P(--)$ exists but it stays below $P_1$ we can consider
 further preimages of the set. We have that $P(-\pm+)= H_+ (P(-\pm))$ will 
be below $P_0$ since we have that $H_+$ maps $S_1$ homeomorphically to 
$\T\times [0,1/2]$. Then this set does not suppose an obstruction to the 
reducibility. On the other hand we have that $P(-\pm)$ define 
two arches around $C_0$ the critical set, then we have
that $P(-\pm-) = H_-(P(-\pm))$ define 
two arches around $P(-)$. It can happen that these arches are 
below $P_1$. If this is the case we can consider again their preimages by $H_-$ and 
$H_+$. The preimages by $H_+$ will be below $P_0$ then they can be discarded, 
because they will not become an obstruction to the reducibility of 
$x_{\alpha,\eps}$. The set $P(-+)$ is below $P_0$, then 
$P(-+-)=H_-(P(-+))$ will be above $P(-)$ and then
 $P(-+--) = H_-(P(-+-))$ will be below $H_-(P(-)) = P(--)$. 
If $P(--)$ does not suppose an obstruction to the reducibility, neither does 
$P(-+--)$. Finally the set $P(----)$ will be an arch above $P_0\cup P(--)$, 
then this can intersect the set $P_1$ becoming an obstruction to the 
reducibility. 

To avoid $P(----)$ being an obstruction to the 
reducibility one should require it to be below $P_1$ for 
any parameter value. This will be an additional constraint 
to the previous ones. In figure
\ref{bounds of reducibility} we have also added this constraint. 

Finally note that the argument used above can be extended to any order. 
Assume that we have $P((-)^{2n})$ which exists but it is not 
an obstruction to reducibility, where $(-)^{2n}$ represents the $2n$ times 
repetition of the symbol $-$. We will have that $H_+(P((-)^{2n}))$ 
can be discarded for being below $P_0$, and $H_+ \circ H_-( P((-)^{2n}))$ 
can be discarded for being below the union of all the set $P((-)^{2k})$ for 
$k\leq n$, then the only set which can suppose an additional restriction 
is $H_- \circ H_-(P((-)^{2n})) = P((-)^{2(n+1)})$. 

Note that the higher order conditions does not necessarily suppose an improvement 
of the previous ones. For example in figure 
\ref{bounds of reducibility} we have considered the additional constraints, but 
not until $P((-)^{12})$ we have had an improvement to the constraint given 
by $P((-)^{4})$.

%  -----------------------------------------------------

%%%%%%%%%%%%%%%%%%%%%%%%%%%%%%%%%%%%%%%%%%%5555
\subsection{Remarks on the constraints} 
%%%%%%%%%%%%%%%%%%%%%%%%%%%%%%%%%%%%%%%%%%%%

%In this subsection we will comment on the estimates 
%obtained before and their consequences. 

%\begin{enumerate} 
%\item 
%First of all let us do two remarks on the discussion done above. 
We have assumed that the sets $P((-)^k)$ for $k$ odd intersect 
the set $P_1$ in exactly $2$ points, but we actually have that this 
assumption can be omitted. If the curves intersect in an even 
number of points we have that the discussion above is still valid 
but considering the different arches at
the same time. When they intersect in an odd number of points 
there is at least one point where both curves are tangent. 
Then the preimage of the intersection  is a single 
point of $P_0$, therefore it can be omitted.
On the other hand, we have only taken into account the case where 
an arch of $P((-)^{2k})$ intersects the set $P_1$
as a possible obstruction to the reducibility
of the curve $x_{\alpha,\eps}$ . It might 
also happen that an arch of the set $P((-)^{2k})$ intersects another
arch of the  set $P((-)^{2r-1})$ for some $k> r > 0$. This can produce also an 
obstruction to the reducibility of the invariant curve $x_{\alpha,\eps}$ as well.  
But this case can be omitted because,  if it occurs, then we can consider 
the ($2r +2$)-th image of both sets by the FLM an we will have that 
$P((-)^{2(k-r-1)}) \cap P_1 \neq \emptyset$.

Some other properties can be deduced for the pre-critical sets. For 
example, if $P((-)^{k_0})=\emptyset$ for some $k_0$, then 
we have that $P((-)^{k})= \emptyset $ for any $k\geq k_0$. On the other hand 
we have that if $P((-)^{2k_0}) \cap P_1 \neq \emptyset$ for some $k_0$ then 
$P((-)^{k})\neq \emptyset$ for any $k\geq k_0$.  Another interesting property is
 the fact that when the set $P((-)^{2k})$ does not intersect 
 $P_1$ for any $k >0 $ then we have that between them there 
exist an invariant compact subset in $\T\times[0,1]$. The set
is delimited on the top by the union of the set
$P_1$ and the sets  $P((-)^{2k +1})$ for any $k>0$, and below by 
the unions of the set $P_0$ and the sets $P((-)^{2r})$ for any $r>0$. 
Note that any curve inside these regions will be reducible. Indeed, 
when this invariant subset exists we will say that the FLM has an
{\bf invariant set of reducibility}.

In the previous subsection the constrains on the reducibility have 
been given through geometric conditions on the pre-critical set. 
The numerical computation of the constrains has been done as 
follows. We have fixed a value for $\alpha$, for instance 
$\alpha=3.2$. Then we have localized a value of $\eps$ such that 
$P((-)^{2k})$ exists. We know that for $\eps= 1-1/\alpha$ this is satisfied, but 
for a better numerical stability of the computation 
it has been convenient to consider lower values of $\eps$. 
Once the parameters have been fixed we have looked for the points 
$P((-)^{2k})\cap P_1$, then these points have been continued in $\eps$ until a tangency 
between the curves has been obtained. Then the tangency has been continued in the 
parameter space as a constraint on the reducibility. A priori, the constraints obtained 
do not have to be better than the previous ones, therefore only the  improvements of 
the constraint have been kept. 

\begin{figure}[t]
\begin{center}
\includegraphics[width=12cm]{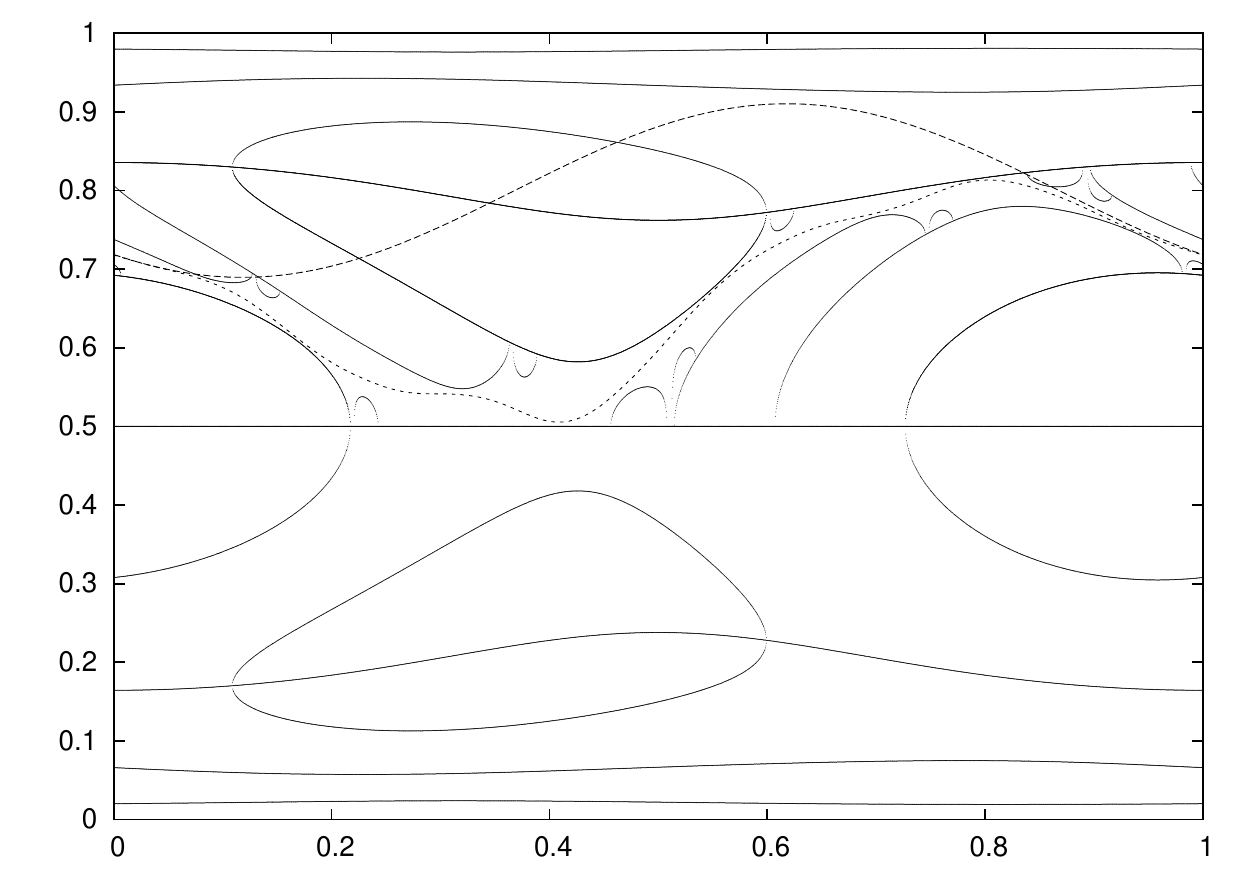}
\caption{ Some remarkable sets of the FLM for $(\alpha,\eps)=(3.2, 0.138)$. 
In a solid line we have the union of the sets $P_0$, $P_{-1}$, $P_{-2}$, $P_{-3}$
and all the sets of the type $P((-)^k)$ for any $k>0$ (recall that we 
admit $P((-)^k)=\emptyset$. In a 
dashed line we have the post-critical set $P_1$ and in a 
dotted line the attracting curve of the map.
The horizontal axis corresponds to $\theta$ and the vertical one to $x$.} 
\label{precritical sets and the attractor} 
\end{center}
\end{figure}

Analyzing the results on the figure \ref{bounds of reducibility}, one can observe
that for values of $\alpha$ close to $2$, the first constraint  
(\ref{A first bound for the reducibility})  coincides with 
the boundary of the parameter space where the 
attractor is a reducible curve. When the parameter $\alpha$ gets further from $2$ the 
constraint is replaced by the constraint given by the set $P((-)^2)$, later on the 
constraint given by the set $P((-)^4)$ becomes optimal, and even further this
one is replaced by the constraint associated to the set $P((-)^{12})$. 
In other words it seems that the minimum of all the possible constrains 
gives the optimal constrain on the reducibility of $x_{\alpha,\eps}$. 
We think this can be due to the fact that whenever there exists an 
invariant set of reducibility, then the 
invariant curve $x_{\alpha,\eps}$ remains in the set. 
We have only found 
numerical evidences for this fact. In figure 
\ref{precritical sets and the attractor} we have plot the 
attractor of the FLM together with  
pre-critical sets of the map. We can observe that, even though the invariant set 
of reducibility is very thin, the invariant curve lays 
completely in its interior. 

Then the different constraints on the reducibility 
would be optimal because they give the region (in the parameter 
space) of existence of the invariant set of reducibility. 
For parameters close to $2$ the invariant set of reducibility 
is delimited only by $P_0$ and $P_1$. 
When the parameter $\alpha$ gets further from $2$, at some point in the 
parameter space the invariant set of reducibility  gets  
delimited by $P_0\cup P((-)^2)$ and $P_1 \cup P(-)$. Then, from that 
point on, the constraint which determines the existence of the set 
is the tangency between $P((-)^2)$ and $P_1$, what 
makes the minimum of both constraints the optimal one for the reducibility. 
When the parameter $\alpha$ gets 
even further from $2$ the existence of the invariant set of reducibility is 
determined by higher order pre-critical sets, and the one which determines the 
existence of the set gives the optimal constrain.

\begin{figure}[t]
\begin{center}
\includegraphics[width=7.5cm]{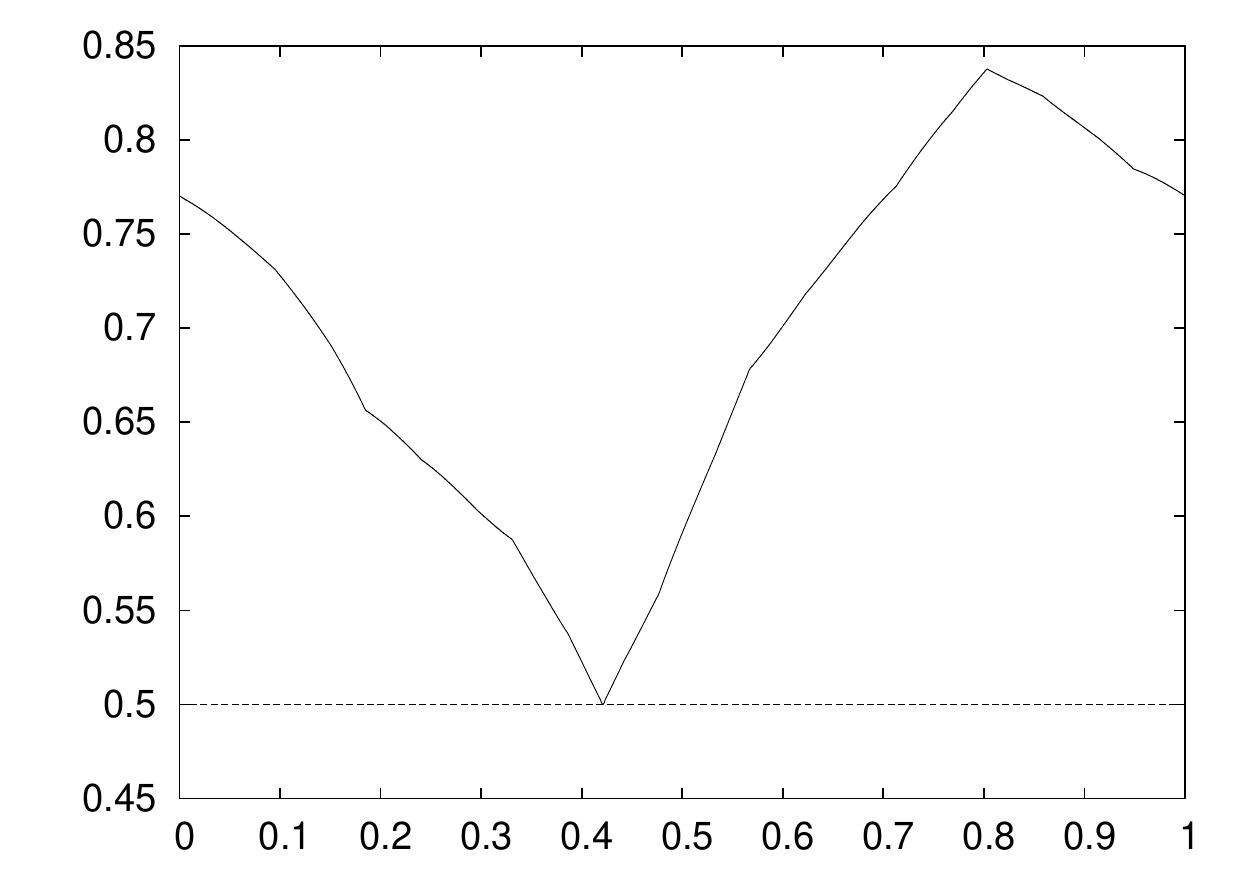}
\includegraphics[width=7.5cm]{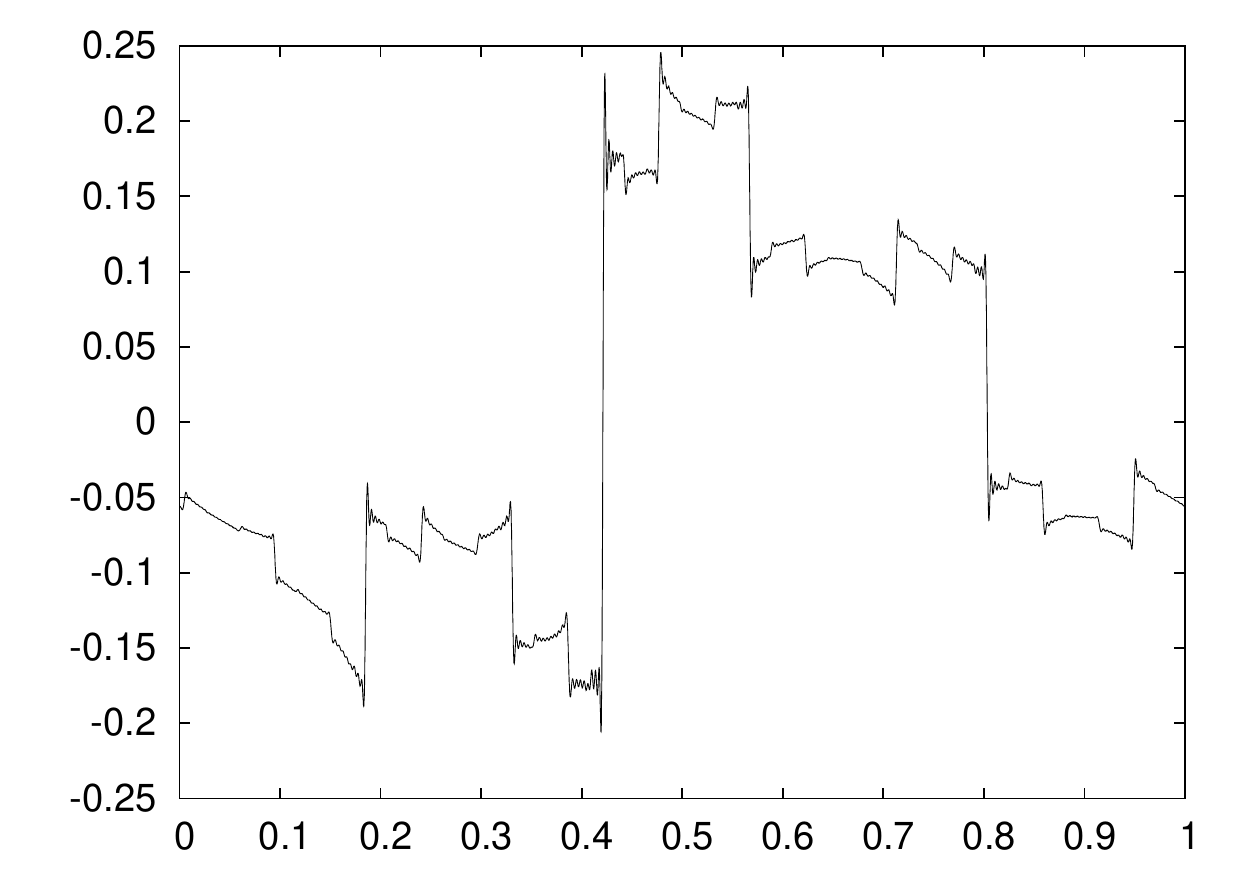}
\includegraphics[width=7.5cm]{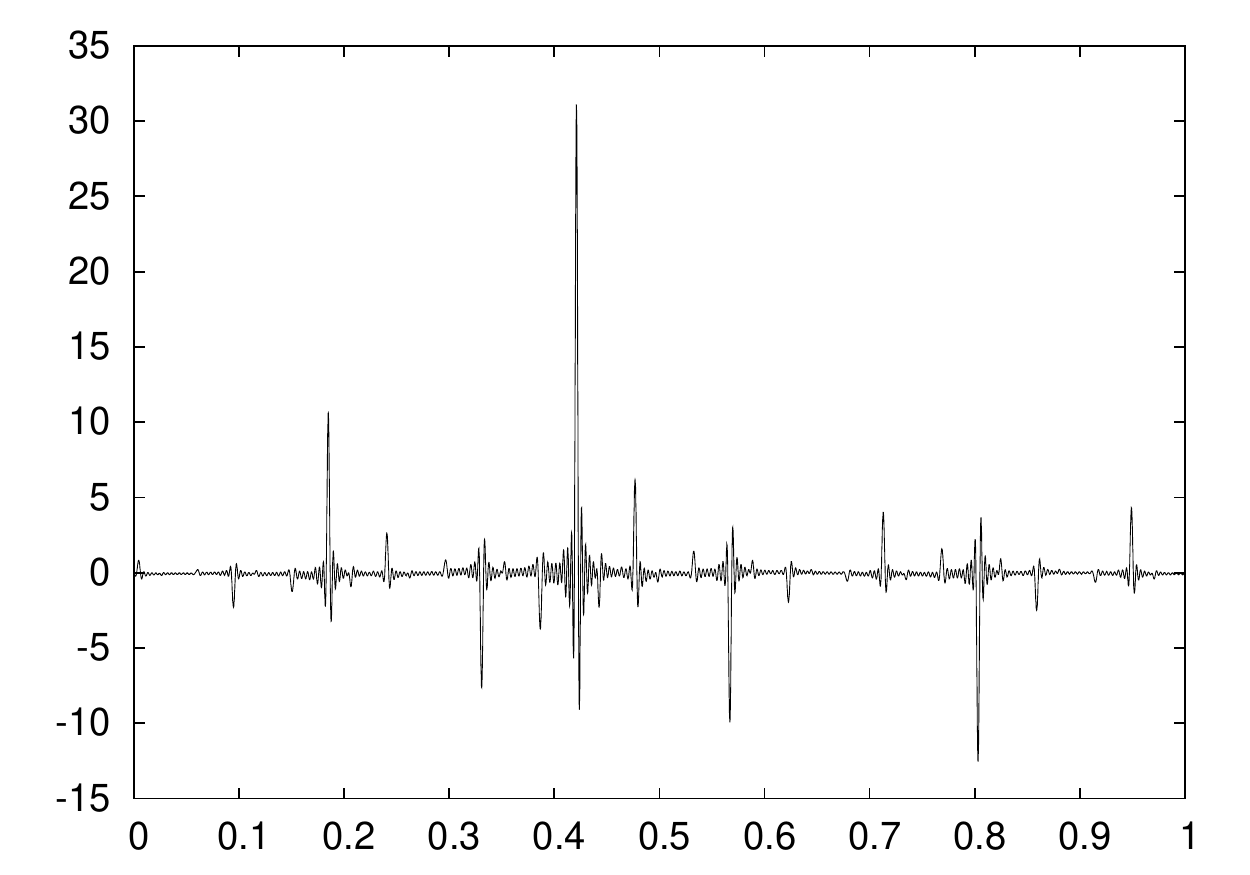}
\includegraphics[width=7.5cm]{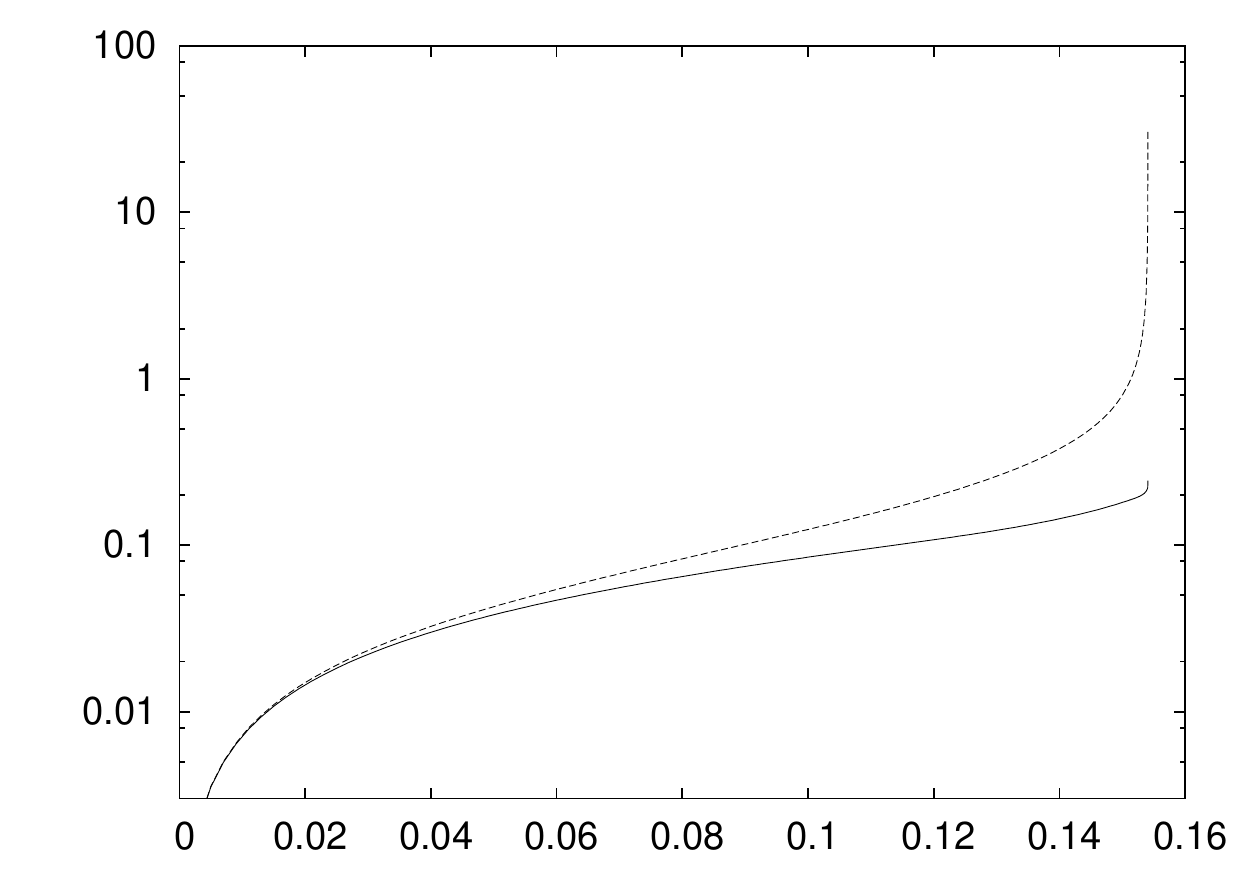}
\caption{ In the top left we have an invariant repelling curve $x(\theta)$ 
of the FLM (\ref{FLM}) for $\alpha$ fixed at $3.5$ and 
$\eps$ just before losing its reducibility (which happens for 
$\eps \approx 0.154086$). In the top right 
and bottom left we have the first ($x'(\theta)$) and the 
second derivative ($x''(\theta)$)  of the same invariant curve.
In these pictures the horizontal axis corresponds to $\theta$. 
In the bottom right we have (in the vertical axis on a logarithmic scale) the 
supremum $|x_\eps'(\theta)|$ (solid line) and $|x_\eps''(\theta)|$ 
(dashed line) as a graph of the parameter $\eps$ (horizontal axis),
where $x_\eps$ is the continuation of the  invariant curve 
$x_0(\theta) = 1-1/\alpha$ for $\alpha$ fixed at $3.5$. }
\label{unstable critical curve}
\end{center}
\end{figure}

%\item 

A final remark is that the constrains obtained 
do not depend on the stability of the invariant curve. 
In other words, they are also valid after the period 
doubling of the attracting set. 
In the unstable case, it is known that the IFT can not be applied 
to prove the persistence of an invariant curve after the reducibility 
is lost (see section 3.5 of \cite{JT05}). Indeed we have 
done numerical computations which suggest that when the invariant curve loses 
its reducibility it suffers a ``fractalization'' process. 
%pproximes the the reducibility-loss the curve is destroyed. 
Concretely, we have continued numerically with respect to $\eps$ the curve 
$x_{\alpha,\eps}$ for $\alpha=3.5$. 
To compute the invariant curve, we have 
approximated it by its truncated Fourier series as described in section 
\ref{section parameter space and reducibility}. 
Looking at the 
bifurcation diagram in figure \ref{FLM parameter space} 
we have that this lays completely to 
the right of the period doubling bifurcation. In figure 
\ref{unstable critical curve} we have displayed the invariant 
curve and its derivative for the last parameter $\eps$ of the 
continuation. We can see that the curve looks quite sharp, due to
big variations of its derivative. In the same figure we have 
$\displaystyle \sup_{\theta\in\T} |x_{\alpha,\eps}'(\theta)|$ and 
$\displaystyle \sup_{\theta\in\T} |x_{\alpha,\eps}''(\theta)|$ against $\eps$. We can 
observe that  
that $\displaystyle \sup_{\theta\in\T} |x_{\alpha,\eps}''(\theta)|$ seems to 
grow unbounded when the parameters of the map get closer to the 
boundary of reducibility. Note that the process of destruction is 
the same as the one obtained  when moving the parameters through the period doubling 
bifurcation. If an unstable curve 
is really destroyed when it loses its reducibility, then we would
have that the constraints on the reducibility  actually enclose the region of 
existence of the curve $x_{\alpha,\eps}$ in the unstable case.

% -----------doubled omega  ----------------------

%\end{enumerate} 

%%%%%%%%%%%%%%%%%%%%%%%%%%%%%%%%%%%%%%%%%%%%%%%
\section{Period doubling and reducibility} 
\label{section period doubling and reducibility} 
%%%%%%%%%%%%%%%%%%%%%%%%%%%%%%%%%%%%%%%%%%%%%%%

In section \ref{section parameter space and reducibility} we have 
studied numerically the bifurcation of the FLM. 
Concretely 
we have  some bifurcation diagrams of the map  
in figure \ref{FLM parameter space}.  
It can be observed that from each parameter value $(\alpha,\eps) = (f_i,0)$
it is born a curve of period doubling bifurcation, which is confined in
a reducibility zone. Moreover, the period doubling bifurcation curve
collides (in a tangent way) with the boundaries of the reducibility, giving 
place to a codimension two bifurcation, namely  
``period doubling - reducibility loss'' bifurcation. 
In the first part of this section we will assume that we have an 
invariant curve and we give a result which says that this tangent 
collision should be expected. 
In the second part we will give a model map for the reducibility regions 
observed in the bifurcations diagram of the FLM. 

%%%%%%%%%%%%%%%%%%%%%%%%%%%%%%%%%%%%%%%%%%%%%%%%%%%%
\subsection{Interaction between  period doubling and reducibility loss 
bifurcations}
%%%%%%%%%%%%%%%%%%%%%%%%%%%%%%%%%%%%%%%%%%%%%%%%%%%%

Theorem 3.3 of \cite{JT05} describes the behavior of 
the Lyapunov exponent $\Lambda = \Lambda(\mu)$ of 
a one parametric family of maps, with respect to 
its parameter $\mu$ close to a value $\mu_0$ 
where the reducibility is lost. Now we will 
present a result which have some  resemblances 
with the cited one, but in the following result 
we consider a two parametric family of maps, and then 
we study the interaction of the period doubling and reducibility loss
bifurcations in the parameter plane.

\begin{thm}
\label{theorem Lyapunov exponent and reducibility}
Consider a one dimensional linear skew product 
\begin{equation}
\label{2-d family of skew products}
\left.
\begin{array}{rcl}
\bar{\theta}&= & \theta  + \omega, \\
\bar{x}& = & a(\theta, \mu,\lambda) x,
\end{array}
\right\}
\end{equation}
which depend ($C^\infty$-)smoothly on $\theta$ and (also ($C^\infty$-)smoothly)
on two parameters $(\mu,\lambda)\in U$ with $U \subset \R^2$ an open set. 

Suppose that there exists a regular curve in the parameter space, parameterized 
by $(\mu_0(t),\lambda_0(t))$ for $t\in [0,\eps_0)$, such that 
the skew product undergoes a loss of reducibility when we cross 
transversally that curve. Let $\Lambda(\mu,\lambda)$ denote  the Lyapunov 
exponent of the invariant curve. Suppose also that 
$\Lambda(\mu_0(0),\lambda_0(0))=0$ and 
$\frac{d}{dt} \Lambda(\mu_0(t)),\lambda_0(t))_{t=0}\neq 0$.

Then there exists a curve in the reducibility zone 
of the parameter space, 
parameterized by $(\mu_1(t),\lambda_1(t))$ for $t\in[0,\eps_1)$,  
with a quadratic tangency with the previous curve and such that 
\[\Lambda(\mu_1(t),\lambda_1(t))=0 \text{ for any } t\in [0,\eps_1).\]
\end{thm}

\begin{proof}
%[Proof of theorem \ref{theorem Lyapunov exponent and reducibility}]
Consider a change of variables to new parameters
$(\bar{\mu}, \bar{\lambda})$ such that the curve $(\mu_0(t),\lambda_0(t))$ goes
to the positive $\mu$ semiaxis. Then the reducibility loss occurs
when $\bar{\lambda}$ crosses transversally zero, for any value for 
$\bar{\mu}\in[0,\bar{\eps})$. We can apply the normal
form near a reducibility loss (Lemma 3.5 of \cite{JT05}) 
for any $\bar{\mu}$ when $\bar{\lambda} =0$.
Then  we have that there exist a change of variables such that 
the system (\ref{2-d family of skew products}) can be transformed to 
the following one
\begin{equation}
\label{2-d family of skew products normal form}
\left.
\begin{array}{rcl}
\bar{\theta}&= & \theta  + \omega, \\
\bar{y}& = & h(\nu(\bar{\lambda})) (\nu(\bar{\lambda}) 
+ \cos(2\pi(\theta-\theta_0(\bar{\lambda})))) y , 
\end{array}
\right\}
\end{equation}
where $\nu(\bar{\lambda})$ is a smooth function 
of the parameter $\bar{\lambda}$, with 
\begin{eqnarray}
\label{nu normal form 1}
\nu(0)& = & 1,
%\nonumber 
\\ 
\label{nu normal form 2} \rule{0pt}{3ex}
\displaystyle \frac{\partial}{\partial \bar{\lambda}} \nu(0) & \neq & 0 ,
\end{eqnarray}
and $h(\nu)$ is a smooth function which never vanishes and 
$\theta_0(\bar{\lambda})$  a smooth function to $\T$.
%Actually the change of variables that we have applied is the one of  
%proposition 
%but omitting the change on the $\theta$-variable).

Note that when applying Lemma 3.5 of \cite{JT05} above, 
the original map depends smoothly 
on the parameter $\bar{\mu}$. Then the change of variables will also depend on it 
and so will do the function $\nu$, $h$ and $\theta_0$. 
Then we can think on these functions depending on both parameters and write  
$\nu=\nu_{\bar{\mu}}(\bar{\lambda})= \nu(\bar{\mu},\bar{\lambda})$, 
$h(\nu)=h_{\bar{\mu}}(\nu_{\bar{\mu}}(\bar{\lambda}))= 
h(\bar{\mu},\nu(\bar{\mu},\bar{\lambda}))$ and $\theta_0(\bar{\lambda})
=\theta_0(\bar{\mu},\bar{\lambda})$. 

The condition (\ref{nu normal form 2}) must be satisfied for 
any $\bar{\mu}\in[0,\bar{\eps})$, then $\frac{\partial}{\partial\bar{ \lambda} }
\nu(\bar{\mu},0)$ can not vanish. Changing $\bar{\mu}$ by $-\bar{\mu}$ 
if necessary, we can suppose that $\frac{\partial}{\partial\bar{ \lambda} }
\nu(\bar{\mu},0)$ is positive. Therefore conditions 
(\ref{nu normal form 1}) and (\ref{nu normal form 2}) can be rewritten as
\begin{eqnarray}
\label{nu normal form 1 B} 
\nu(\bar{\mu},0)& = & 1 \text{ for any } \bar{\mu}\in[0,\bar{\eps}),
%\nonumber
\\ 
\displaystyle \frac{\partial}{\partial \bar{\lambda} }
\label{nu normal form 2 B}
\nu(\bar{\mu},0) & > & 0  \text{ for any } \bar{\mu}\in[0,\bar{\eps}).
\end{eqnarray}
Note that we are assuming that we have reducibility for $\bar{\lambda}>0$ 
and non-reducibility for $\bar{\lambda}<0$.

Using the formula (\ref{integral ln cos}) we can compute the Lyapunov 
exponent of the system, 
\begin{eqnarray} 
\Lambda(\bar{\mu},\bar{\lambda})  
& = & \displaystyle %\frac{1}{2\pi} \int_0^{2\pi} 
\int_0^1
\ln | h(\bar{\mu}, \nu(\bar{\mu},\bar{\lambda})) ( \nu(\bar{\mu},\bar{\lambda}) 
+ \cos(2\pi(\theta - \theta_0(\bar{\mu},\bar{\lambda})))) | d\theta 
\nonumber \\ 
& = & \rule{0pt}{5ex} \displaystyle 
\int_0^{1} \ln |  h(\bar{\mu}, 
\nu(\bar{\mu},\bar{\lambda})) | d\theta 
+  \int_0^{1} \ln | \nu(\bar{\mu},\bar{\lambda}) 
+ \cos(2\pi(\theta - \theta_0(\bar{\mu},\bar{\lambda}))) | d\theta
\nonumber \\ 
&=&
\rule{0pt}{10ex}
\left\{ \begin{array}{lcr}
\displaystyle 
\ln\left|\frac{ h(\bar{\mu}, \nu(\bar{\mu},\bar{\lambda}))}{2}\right| 
  & \text{ if } &  \bar{\lambda} \leq 0 ,\\
\rule{0pt}{5ex} \displaystyle 
\ln\left|\frac{ h(\bar{\mu}, \nu(\bar{\mu},\bar{\lambda}))}{2}\right|  
+ \operatorname{arccosh}|\nu(\bar{\mu},\bar{\lambda})|
& \text{ if } &   \bar{\lambda} \geq 0.
\end{array} \right.
\label{Lyapunov exponent of the normal form}
\end{eqnarray}

By hypothesis we have that $\frac{d}{dt} \Lambda(\mu_0(t)),
\lambda_0(t))_{t=0}\neq 0$, performing the change of variables 
to the original parameters, we have that 
\[
\frac{\partial}{\partial \bar{\mu}} \Lambda(\bar{\mu},
\bar{\lambda})_{(\bar{\mu},\bar{\lambda})=(0,0)} = \frac{d}{dt} \Lambda(\mu_0(t)),
\lambda_0(t))_{t=0}\neq 0.
\]

If we differentiate equation (\ref{Lyapunov exponent of the normal 
form}) with respect to $\bar{\mu}$ with the constraint $\bar{\lambda}=0$ we have 
\[
\frac{ \partial }{ \partial \bar{\mu} } \Lambda(\bar{\mu},0) = 
\frac{ \frac{\partial}{\partial \bar{\mu}} ( h(\bar{\mu},0)) }
{ h(\bar{\mu},\nu(\bar{\mu},0)) } = 
\frac{ \partial_{\bar{\mu}} h(\bar{\mu},0) + 
\partial_\nu h(\nu(\bar{\mu},0)) \cdot \partial_{\bar{\mu}} \nu(\bar{\mu},0) }
{ h(\bar{\mu},\nu(\bar{\mu},0)) }.
\]
Differentiating in (\ref{nu normal form 1 B}) we have 
$\partial_{\bar{\mu}} \nu(\bar{\mu},0) = 0$, then from the last two equations it follows
that
\begin{equation}
\label{equation auxiliar proof normal forms}
\frac{ \partial_{\bar{\mu}} h(0,\nu(0,0)) }{ h(0,\nu(0,0)) } 
\neq 0.
\end{equation}

Note that for $\bar{\lambda} >0$ we have that 
\[ \Lambda(\bar{\mu},\bar{\lambda})=0 \text{ if, and only if, } 
|\nu(\bar{\mu},\bar{\lambda})|= \cosh
\left(\ln{\frac{|h(\bar{\mu},\nu(\bar{\mu}, \bar{\lambda}))|}{2} }\right). 
\]
As $\nu(\bar{\mu}, 0)= 1$, by continuity we have 
$|\nu(\bar{\mu}, \bar{\lambda})|= \nu(\bar{\mu}, \bar{\lambda})$ 
for any $\bar{\lambda}$ small enough. But 
then, recall the 
$h(\bar{\mu},\nu(\bar{\mu}, \bar{\lambda}))$ is a smooth function which never 
vanishes (in a neighborhood of $\bar{\lambda}=0$). Then, if we consider the function 
\begin{equation}
\label{equation definition of rho}
\rho(\bar{\mu},\bar{\lambda}): = \nu(\bar{\mu},\bar{\lambda}) - \cosh
\left(\ln{\frac{|h(\bar{\mu},\nu(\bar{\mu}, \bar{\lambda}))|}{2}} \right),
\end{equation} 
we have that it is a smooth function on the parameters $(\mu,\lambda)$ for 
$\bar{\mu}\in[0,\bar{\eps})$ and $\bar{\lambda}$ small enough. Moreover 
from (\ref{Lyapunov exponent of the normal form}) we have that
(when $\bar{\lambda}<0$)  $\Lambda(\bar{\mu},\bar{\lambda})=0$ if, and only if, 
$\rho(\bar{\mu},\bar{\lambda})=0$. 

To prove the theorem we first will apply the IFT to the function 
$\rho(\bar{\mu},\bar{\lambda})$, then we will
 check that the curve that we obtain is tangent to the set $\bar{\lambda}=0$ 
and finally that it belongs to the upper semi-plane $\bar{\lambda}>0$. 

From equation  (\ref{Lyapunov exponent of the normal form}) and the 
hypothesis $\Lambda(0,0)=0$ it follows that 
\begin{equation}
\label{equation ln h equal 0}
\ln{\frac{|h(0,\nu(0, 0))|}{2}}=0.
\end{equation}
Therefore using (\ref{nu normal form 1 B}) it follows that $\rho(0,0) = 0$.

To apply the IFT we need to compute first the derivative of $\rho$ with 
respect to $\lambda$ and check that it is different from zero. We have 
\[
\frac{\partial }{\partial \bar{\lambda}} \rho(\bar{\mu},\bar{\lambda}) = 
\partial_{\bar{\lambda}} \nu(\bar{\mu},\bar{\lambda}) - \sinh
\left(\ln{\frac{|h(\bar{\mu},\nu(\bar{\mu}, \bar{\lambda}))|}{2}} \right) 
\frac{ \partial_\nu h(\bar{\mu}, \nu(\bar{\mu},\bar{\lambda})) 
\cdot \partial_{\bar{\lambda}} \nu(\bar{\mu},\bar{\lambda}) }
{ h(\bar{\mu},\nu(\bar{\mu},\bar{\lambda})) }.
\]
Evaluating at $(\bar{\mu},\bar{\lambda})=(0,0)$ and using 
(\ref{nu normal form 2 B}) we have 
\[
\frac{\partial }{\partial \bar{\lambda}} \rho(0,0) = 
\partial_{\bar{\lambda}} \nu(0,0) > 0.
\] 
Then we have that there exist an interval $[0,\bar{\eps}_1)$ and 
a function $\lambda:[0,\bar{\eps}_1) \rightarrow \R$ such that 
$\lambda(0)=0$ and $\Lambda(\mu,\lambda(\mu))=0$ for any 
$\mu\in[0,\bar{\eps}_1)$. Moreover we have that 
\begin{equation}
\label{equation derivative of lambda of mu}
\bar{\lambda}'(\bar{\mu}) = -
\frac{\partial_{\bar{\mu}} \rho(\bar{\mu},\bar{\lambda}(\bar{\mu}))}
{\partial_{\bar{\lambda}} \rho(\bar{\mu},\bar{\lambda}(\bar{\mu}))}.
\end{equation} 

To finish with the proof we only have to check that $\lambda'(0)=0$ and 
$\lambda''(0)>0$. Differentiating in (\ref{equation definition of rho}) with respect 
to $\mu$ we have, 
\begin{equation}
\label{equation partial rho respect mu}
\frac{\partial }{\partial \bar{\mu}} \rho(\bar{\mu},\bar{\lambda}) = 
\partial_{\bar{\mu}} \nu(\bar{\mu},\bar{\lambda}) - \sinh
\left(\ln{\frac{|h(\bar{\mu},\nu(\bar{\mu}, \bar{\lambda}))|}{2}} \right) 
\frac{ \partial_{\bar{\mu}} h(\bar{\mu}, \bar{\lambda}) + 
\partial_\nu h(\bar{\mu}, \nu(\bar{\mu},\bar{\lambda})) 
\cdot \partial_{\bar{\mu}} \nu(\bar{\mu},\bar{\lambda}) }
{ h(\bar{\mu},\nu(\bar{\mu},\bar{\lambda})) }. 
\end{equation}
Using equations (\ref{nu normal form 2 B}) and (\ref{equation ln h equal 0}) 
we can evaluate the last equation 
at $(\bar{\mu},\bar{\lambda}) = (0,0)$ and we 
have $\frac{\partial }{\partial \bar{\mu}} 
\rho(\bar{\mu},\bar{\lambda})=0$, therefore (using 
(\ref{equation derivative of lambda of mu})) we have $\bar{\lambda}'
(\bar{\mu})=0$.

For the condition $\lambda''(0)>0$, differentiating on equation (\ref{equation 
derivative of lambda of mu}) and evaluating at $0$ 
(and using also $\bar{\lambda}'(\bar{\mu})=0$) we have that
\[
\bar{\lambda}''(0) = - 
\frac{\frac{\partial^2}{\partial\mu^2} \rho(0,0)}
{\partial_{\bar{\lambda}} \rho(0,0)} .
\]
Finally, from (\ref{nu normal form 1 B}) we have also that 
$\frac{\partial^2}{\partial\mu^2} \nu(0,0) = 0$, then differentiating 
equation (\ref{equation partial rho respect mu}), evaluating 
at $(\bar{\mu},\bar{\lambda})=(0,0)$ and simplifying all the zero terms we have
\[
\frac{\partial^2}{\partial\mu^2} \rho(0,0) = - 
\left( 
\frac{ \partial_{\bar{\mu}} h(0,\nu(0,0)) }{ h(0,\nu(0,0)) }
\right) ^2,
\]
which is different from zero due to (\ref{equation auxiliar proof normal 
forms}), therefore $\bar{\lambda}''(0)>0$.
\end{proof}

%%%%%%%%%%%%%%%%%%%%%%%%%%%%%%%%%%%%%%%%%%%%%%%%55
\subsection{A model for the reducibility regions}
%%%%%%%%%%%%%%%%%%%%%%%%%%%%%%%%%%%%%%%%%%%%%5%%%%
 
In this section we  present  a two parametric skew product
to model the interaction between the reducibility loss and 
the period doubling bifurcation 
%(and such that it uncouples for one of its parameters). 
%such that it uncouples to a one dimensional map for 
%one of this parameters, it has a period doubling bifurcation 
%of its attracting set and the invariant curve $x=0$ has 
%a reducibility-loss. 
We will see that this map is helpful for the understanding 
of the reducibility regions observed in the parameter space of the FLM. 

Let us recall first the period doubling (and pitchfork) bifurcation 
for one dimensional maps. Consider the one dimensional  map 
$\bar{x}= x(x^2- \mu)$, where $\mu$ is a parameter in the real 
line. It is well known that this map undergoes pitchfork 
bifurcation when the parameter $\mu$ crosses the value $-1$ 
and it undergoes a period doubling bifurcation when $\mu$ 
crosses the value $1$.

The model map that we propose is the following,
\begin{equation}
\label{model q.p period doubling}
\left.
\begin{array}{rcl}
\bar{\theta} & = &\theta + \omega ,\\
\bar{x} & =  & x\left(x^2 - (\mu + \lambda \cos(2\pi\theta))\right) ,
\end{array}
\right\}
\end{equation}
where $\mu$ and $\lambda$ are parameters, and $\omega$ is Diophantine.

Note that when $\lambda$ is equal to zero the system
uncouples and we obtain the model for the generic unfolding of the period 
doubling bifurcation when $\mu$ crosses $1$. Indeed the q.p. forcing 
has been considered in such a way that the set $x=0$, namely the 
trivial invariant set, is always  an invariant curve of the map 
(for any parameter values). 

The linear dynamics around this trivial invariant set is
\begin{equation}
\label{skew trivial inviariant set}
\left.
\begin{array}{rcl}
\bar{\theta} & = &\theta + \omega ,\\
\bar{x} & =  &  - (\mu + \lambda \cos(2\pi\theta)) . 
\end{array}
\right\}
\end{equation}
We are in the $C^\infty$ framework, therefore by Corollary 1 
in \cite{JT05} we have that the trivial invariant set 
is reducible if, and only if, $|\lambda| < |\mu|$. In other words, the 
loss of reducibility bifurcation correspond to the lines in 
the parameter space given by $\lambda= \pm \mu$.

On the other hand,  we can compute $\Lambda=\Lambda(\mu,\lambda) $
the Lyapunov exponent of the trivial invariant curve $x=0$.  
Recall that, 
\begin{equation}
\label{integral ln cos}
%\frac{1}{2\pi} \int_0^{2\pi} 
\int_0^1 \ln | \tau + \cos(2\pi\theta) | d\theta = 
\left\{ \begin{array}{lcr}
-\ln{2}  & \text{ if } &  |\tau | \leq 1,\\
-\ln{2}  + \operatorname{arccosh} |\tau| & \text{ if } &  |\tau | \geq 1,
\end{array} \right.
\end{equation}
and $\displaystyle 
\operatorname{arccosh} |\tau| = \ln\left(|\tau|+\sqrt{\tau^2-1}\right)$ for 
$|\tau|\geq 1$. Using this it is easy to check that
\[
\begin{array}{rcl} 
\Lambda(\mu,\lambda) & = & \displaystyle %\frac{1}{2\pi} \int_0^{2\pi} 
\int_0^1 \ln | \mu + \lambda \cos(2\pi\theta)) | d\theta \\ & = &
\rule{0pt}{10ex}
\left\{ \begin{array}{lcr}
\displaystyle 
\ln{\left(\frac{|\mu| + \sqrt{\mu^2-\lambda^2}}{2}\right)}
  & \text{ if } &  |\lambda| \leq |\mu|,\\
\rule{0pt}{5ex} \displaystyle 
-\ln{\left|\frac{\lambda}{2}\right|}  & \text{ if } &   |\lambda| \geq |\mu|,
\end{array} \right.\end{array}
\]

Concretely, 
for $|\lambda| \leq |\mu|$ we have that the skew product associated to 
the trivial invariant set $x= 0$ is reducible to 
\[
\left.
\begin{array}{rcl}
\bar{\theta} &= & \theta  + \omega, \\
 \rule{0pt}{3ex} 
\displaystyle \bar{x} &= & \displaystyle - 
\frac{\mu + \operatorname{sign}(\mu) \sqrt{\mu^2-\lambda^2}}{2} x.  \\
\end{array}
\right\}
\] 
Therefore, in the reducible case ($|\lambda| \leq |\mu|$) 
 the trivial invariant set changes its stability when 
\[ |\mu| = 1 + \frac{\lambda^2}{4}. \]

One should expect that the map (\ref{model q.p period doubling})
 has a pitchfork bifurcation when one crosses the parabola 
$\mu= - 1 -\frac{\lambda^2}{4}$ and respectively 
a period doubling bifurcation when crosses $\mu = 1 + \frac{\lambda^2}{4}$. 
Unfortunately this can only be proved in a small region 
of the reducibility zone. More concretely we have the following 
results. 

\begin{prop}
\label{proposition existense of period doubling}
The model map (\ref{model q.p period doubling}) has
\begin{itemize} 
\item Two-periodic  (continuous and non-trivial) invariant curves when $0<|\lambda|<\mu$ and
$1+\frac{\lambda^2}{2}< \mu  < \frac{3}{2} - 2|\lambda|$.
\item Two (continuous and non-trivial) invariant curves when 
$\mu<-|\lambda|<0$ and $\frac{-3}{2} + 2|\lambda| < 
\mu <- 1+\frac{\lambda^2}{2} $.
\end{itemize}
\end{prop}

\begin{figure}[t]
\begin{center}
\includegraphics[width=10cm]{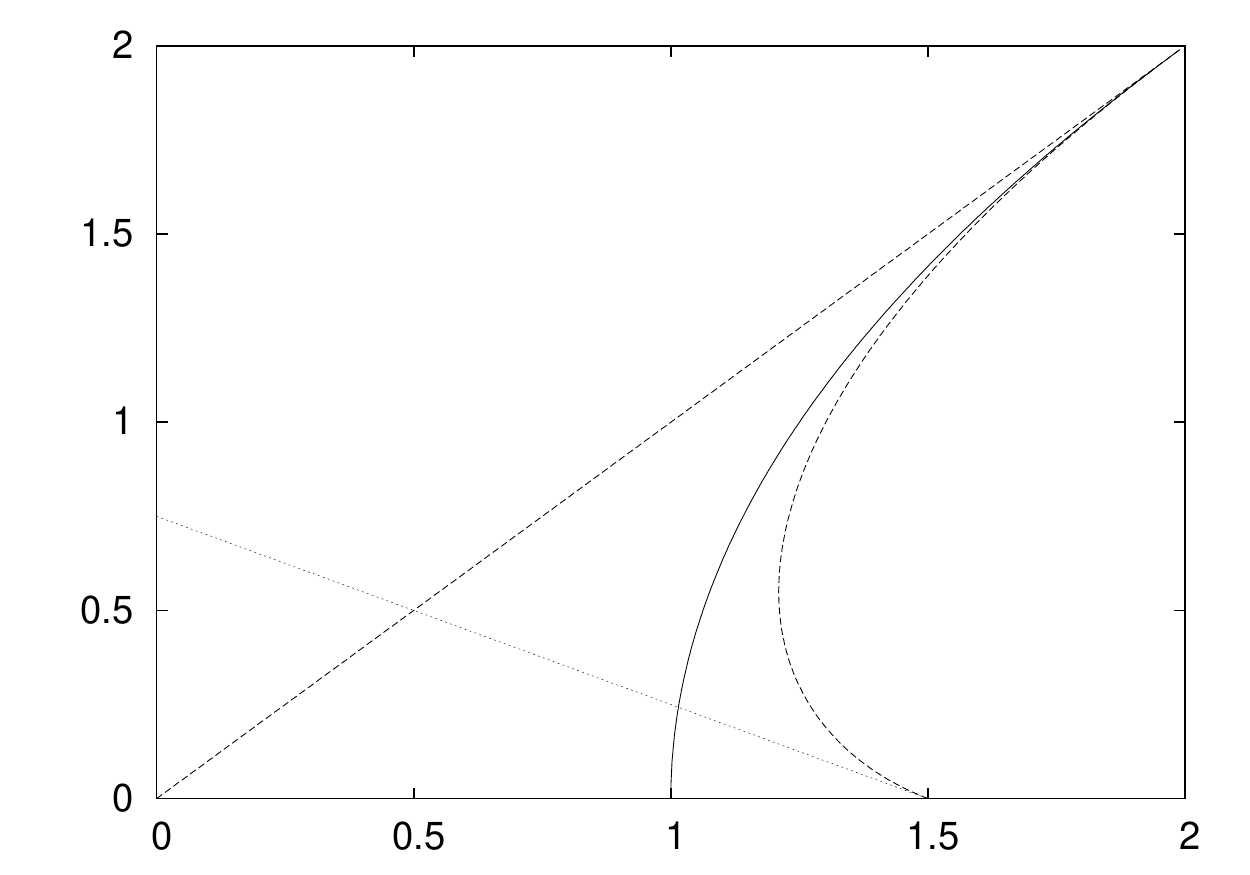}
\caption{Different remarkable curves of the map 
(\ref{model q.p period doubling}) are displayed. In a solid line 
we have the change of stability of the trivial invariant set. In a 
dashed line the loss of reducibility of the attracting set, 
which is the trivial invariant set or a two periodic 
invariant curve. Finally we have in a dotted line the boundary 
of applicability of the proposition 
\ref{proposition existense of period doubling}.
The horizontal axis corresponds to the parameter $\mu$ and 
the vertical one to $\lambda$.  
}
\label{period doubling model parameter space}
\end{center}
\end{figure}

\begin{proof}
%[Proof of proposition \ref{proposition existense of period doubling}]
%The proof will follow using the theorem \ref{theorem Jager existence invariant curves}. 
%
We will prove only the first item of the proposition, since the other 
one is completely analogous. 

We  are interested  on the two periodic invariant curves 
of (\ref{model q.p period doubling}). 
Let us denote by $f$ the function which defines the system, i.e. 
 $f(\theta,x) = x\left(x^2 - (\mu + \lambda \cos(2\pi\theta))\right)$. 
Observe that $f(\theta,-x) = - f(\theta,x)$ for any $(\theta,x)\in 
\T \times \R$ and any parameters values. It follows that, 
if $u:\T \rightarrow \R$ is an invariant curve of (\ref{model q.p period doubling}), 
then $-u$ also is. Moreover, given a curve $u$ (different 
from the trivial invariant set $u\equiv 0$)  satisfying
\[u(\theta+\omega)=-f(\theta, u(\theta)), \text{ for any } \theta \in \T, \] 
we have that  $u$ is a two periodic solution.

Let us consider the system 
\begin{equation}
\label{model q.p period doubling B}
\left.
\begin{array}{rcl}
\bar{\theta} & = &\theta + \omega \quad ,\\
\bar{x} & =  & h(\theta,x) ,
\end{array}
\right\}
\end{equation}
where $h(\theta,x) = -f(\theta,x)= - 
x\left(x^2 - (\mu + \lambda \cos(2\pi\theta))\right)$. It is 
clear that an invariant curve (different from $x(\theta)\equiv 0$)
for the system above is a two periodic solution of 
(\ref{model q.p period doubling}). We will apply  theorem 4.2 of 
\cite{Jag03}  
to this map in order 
to prove the existence of periodic solutions of (\ref{model q.p period doubling}). 
Let us check the conditions of the theorem. We have that 
\[ \frac{\partial}{\partial x} h (\theta,x) = 
-3 x^2 + \mu + \lambda \cos(2\pi\theta). \] 
When $0\leq |\lambda| < \mu$ we have that 
$\frac{\partial}{\partial x} h (\theta,x) >0 $ for any $(\theta, x) $ 
such that $x^2 < \frac{\mu + \lambda \cos(2\pi\theta)}{3}$. Concretely 
this condition is satisfied when $x^2 < \frac{\mu - |\lambda|}{3}$. 

On the other hand, we have that 
\[  \frac{\partial^2}{\partial x^2} h (\theta,x) =  -6 x \text{ and }
\frac{\partial^3}{\partial x^3} h (\theta,x) = -6, \] 
therefore 
\[ S_x h(\theta,x) = 
\frac{-6} {\frac{\partial}{\partial x} h (\theta,x)} - 
\frac{3}{2}\left( \frac{-6 x } 
{\frac{\partial}{\partial x} h (\theta,x)}\right)^2. \]
Note that whenever
$ \frac{\partial^2}{\partial x^2} h (\theta,x) >0 $ then we have that 
$S_x h(\theta,x) <0$.  Finally the only hypothesis remaining to be
satisfied is to find an interval $I$ such that $\T\times I$ is 
invariant by the system (\ref{model q.p period doubling B}). 
Let us consider the interval 
$I=\left[-\sqrt{\frac{\mu - |\lambda|}{3}} + \delta,\sqrt{\frac{\mu - |\lambda|}{3}}- \delta
\right] $, with $\delta$ an arbitrary small value yet to be defined. 
 Note  that for any point on the interval 
the monotonicity condition is satisfied, then it is enough to check that
the interval satisfies the invariance condition. 
In order to have invariance of the interval we have to check that
\[
\left| h\left(\theta,x \right)\right| < 
\sqrt{\frac{\mu - |\lambda|}{3}},
\]
for any $(\theta,x) \in \T \times I$. Since the monotonicity condition 
is satisfied, it is enough to check that 
\[
\left| h\left(\theta, \pm \sqrt{\frac{\mu - |\lambda|}{3}} -\delta\right)\right| < 
\sqrt{\frac{\mu - |\lambda|}{3}} - \delta.
\]
Then it follows that it is enough to have 
\[\left|\frac{2\mu}{3} - \left(\frac{|\lambda|}{3} + \lambda \cos(2\pi\theta)\right) - \delta
\right| < 1. \]
when $\mu + 2|\lambda| < \frac{3}{2}$, we have that there exists a value 
of $\delta$ sufficiently small such that this is satisfied for any $\theta\in\T$. 

We are in situation of applying the theorem 4.2 of \cite{Jag03}, 
 we know 
that $u(\theta) \equiv 0$ is always a continuous invariant curve, which 
is contained in the set $\T\times I$. Moreover we know its 
Lyapunov exponent explicitly, therefore when this crosses zero, 
the theorem implies that there exist two invariant curves (with 
negative Lyapunov exponents) of the system (\ref{model q.p period 
doubling}), which correspond to periodic solutions of the system 
(\ref{model q.p period doubling B}). Moreover the invariant curves  
have negative Lyapunov exponent, therefore the periodic invariant 
curve of the original map is attracting. 
\end{proof}

To illustrate this last 
proposition, in figure \ref{period doubling model parameter space} 
we have plotted the curve which constrains the validity of proposition 
\ref{proposition existense of period doubling}. 

Let us assume that in the reducible case ($|\lambda|<|\mu|$) the parabola 
$\mu= 1+ \frac{\lambda^2}{4}$ corresponds to a period doubling 
bifurcation. Note that this parabola has a tangency at the points 
$(\mu,\lambda) = (2, \pm 2)$ with the boundary of reducibility,
as predicted by theorem \ref{theorem Lyapunov exponent 
and reducibility}.  Then these points 
in the parameter space would correspond to the 
``period doubling - reducibility loss'' bifurcation, since it is the 
point where both curves  merge. 

Recall that in section \ref{section parameter space and reducibility} 
we have reported how  the period doubling bifurcation of the FLM were 
enclosed inside regions of reducibility. In the proposed model 
(\ref{skew trivial inviariant set}) for the reducibility regions 
it is easy to justify this behavior.

Differentiating the equation which defines the map 
we have that the critical region of the map is 
$P_0 = \{ (\theta,x) \in \T\times I  | \thinspace 3 x^2 - 
(\mu + \lambda \cos(2\pi\theta)) =0 \} $. Following the arguments of section 
\ref{section obstruction to reducibility} we have that a reducible 
invariant (or periodic) curve of  (\ref{skew trivial inviariant set}) 
can not intersect the critical set. Note that for $|\lambda|<|\mu|$ the 
set  $P_0$ is composed by two closed curves in the cylinder,  one in each side 
of the trivial invariant set. Moreover when  $|\lambda|$ tends to $ |\mu|$ we 
have that the components get closer to the trivial invariant set $x=0$. 

If we assume that the parabola in the parameter space $\mu= 1+ \frac{\lambda^2}{4}$  
corresponds  to a period doubling bifurcation, 
then we have that when the trivial invariant 
set becomes unstable then a period two solution must be created in 
its neighborhood. But then, when $|\lambda|$ tends to $|\mu|$ we 
have that the set $C_0$ gets closer and closer to the trivial
 invariant set, therefore there is no room for the period doubled 
curve to be reducible. This explains why arbitrarily close to 
the  ``period doubling - reducibility loss'' bifurcation parameter 
$(\mu,\lambda) = (2, \pm 2)$ one can observe a reducibility loss of the 
period two invariant curve. 

In figure \ref{period doubling model parameter space} 
there are shown the different bifurcation curves of the map 
(\ref{skew trivial inviariant set}). The reducibility 
of the period two invariant curve have been estimated numerically. 
Let us remark the resemblance of the region of reducibility of the 
attractor with the same regions of the bifurcation diagram of the FLM.

%%%%%%%%%%%%%%%%%%%%%%%%%%%%%%%%%%%%%%%%%%%%%%%%%%%%%%%%%%%%%%%%%%%%%%%%%
%%%%%%%%%%%%%%%%%%%%%%%%%%%%%%%%%%%%%5
\section{Summary and conclusions}
\label{section summary and conclusions}
%%%%%%%%%%%%%%%%%%%%%%%%%%%%%%%%%%%%%%%

We study the Forced Logistic Map as a toy model for
the truncation of the period doubling cascade of invariant curves.
In section \ref{section parameter space and reducibility}
we have done a numerical analysis of the bifurcation diagram of the FLM,
which is displayed in figure \ref{FLM parameter space}. This computation
revealed that each period doubling bifurcation curve in the
parameter space is confined inside a region
where the attracting invariant curve is reducible.
Now we can use studies done
in sections \ref{section obstruction to reducibility}  and 
\ref{section period doubling and reducibility}
to review the analysis of the bifurcation diagram 
of the FLM done in section \ref{section parameter space and reducibility}.

In section \ref{section obstruction to reducibility} we have done
a study of the critical set, their images and their preimages.
We have constructed different constrains
in the parameter space for the
reducibility of the invariant curve $x_{\alpha,\eps}$ (which is
the continuation of the invariant curve
$x_{\alpha,0}(\theta)= 1-1/\alpha$ for $\eps>0$). We have
 also illustrated how the combination of all these constrains seems to
be the optimal constrain for the reducibility of the invariant curve.
In other words, they approximate the boundary of reducibility
of the attracting set of the map. Using the notation introduced in
section  \ref{section parameter space and reducibility}, we have
that these constraints characterize the curve $C^+_{0}$.
In the bifurcation diagram of the figure \ref{FLM parameter space}
only the properties
of the stable set are reflected, but we have that the constrains are still
valid after the period doubling. Actually we conjecture that these constrains 
give the boundary of existence of the  curve $x_{\alpha,\eps}$ when
it is unstable.

In section \ref{section period doubling and reducibility} we have
studied the interaction between the reducibility loss and the 
period doubling bifurcation. For the case of
linear skew products we have theorem \ref{theorem Lyapunov 
exponent and reducibility}, which says that generically
we can expect the period doubling bifurcation curves and the
reducibility loss bifurcation curve to be tangent.
The diagram of the figure \ref{FLM parameter space}
has been done in terms of the attracting set of
the FLM. As the FLM is not a linear skew product, this theorem
is not applicable. But it can be applied to the linear skew
product given by the linearization of the map around the invariant curve.
In the same section, we have also given a model for the interaction of
the period doubling bifurcation and the reducibility loss. With this model
we have seen that, if the period doubling is close to a reducibility loss,
then there is an obstruction to the reducibility of the period doubled curve. 
This explains why the curves $C^+_{0}$, $D_1$ and $C^-_{1}$ meet at the same point
(using again the notation of section \ref{section parameter space and 
reducibility}). Moreover,
theorem \ref{theorem Lyapunov exponent and reducibility}
gives a good explanation of why they do it
in a tangent way. Finally, let us remark that the study done in section
\ref{section period doubling and reducibility}
does not depend on the map considered, therefore
it can be extended to the rest of reducibility regions
determined by $C^+_{i-1}$ and $C^-_{i}$  (and containing $D_i$).

Part of the study done here will 
be continued in \cite{JRT11a, JRT11b, JRT11c}. 
In \cite{JRT11a} we will propose an  
extension of the renormalization the theory 
for the case of one dimensional quasi-periodic forced maps. Using this 
theory we will be able to prove that the curves $C^{\pm}_i$ of 
reducibility loss bifurcation really exists (for $\eps$ small enough). 
In \cite{JRT11b} we will use the theory proposed in the previous one 
to study the asymptotic behavior of the reducibility loss bifurcations 
when the period goes to infinity. In the previous two articles several 
conjectures will be done. In \cite{JRT11c} we will support numerically this 
conjectures. We will give also numerical 
evidences of the self-renormalizable character of the bifurcation 
diagram of the figure \ref{FLM parameter space} when different values 
of the rotation number $\omega$ are considered.

%\bibliography{bibbase.bib}
\bibliographystyle{plain}

\end{document}